\documentclass[11pt]{article}
\usepackage{latexsym,amssymb,amsmath,amsfonts,amsthm}
\usepackage{graphicx}
\usepackage{subfigure}
\usepackage{placeins}
\usepackage{color}
\usepackage{indentfirst}
\usepackage{mathrsfs}
\usepackage{amsmath,mathrsfs,amsfonts}
\usepackage{tikz}
\usepackage{lineno}
\usepackage{xcolor}
\usepackage{standalone}
\usepackage{cite}
\numberwithin{equation}{section}
\newcommand{\no}{\nonumber}
\newcommand{\be}{\begin{eqnarray}}
	\newcommand{\ben}{\begin{eqnarray*}}
		\newcommand{\en}{\end{eqnarray}}
	\newcommand{\enn}{\end{eqnarray*}}

\newtheorem{theorem}{Theorem}[section]
\newtheorem{lemma}[theorem]{Lemma}
\newtheorem{corollary}[theorem]{Corollary}

\newtheorem{remark}[theorem]{Remark}

\definecolor{rot}{rgb}{0.000,0.000,0.000}
\definecolor{blau}{rgb}{0,0,1}

\begin{document}
	\renewcommand{\theequation}{\arabic{section}.\arabic{equation}}
	\begin{titlepage}
		
		\title{\bf {\color{rot}{Well-posedness of grating diffraction problems for plane wave incidence: explicit dependence on wavenumbers and incident angles}}}

		\author{Linlin Zhu \footnote{Schoole of Mathematics Sciences and LPMC, Nankai University, Tianjin 300071, People's Republic of China (llzhu@mail.nankai.edu.cn).} ,\ \ Guanghui Hu \footnote{(Corresponding author) School of Mathnematical Sciences and LPMC, Nankai  University, Tianjin 300071, People's Republic of China (ghhu@nankai.edu.cn).}}

		\date{}
	\end{titlepage}
	\maketitle
	\vspace{.2in}
	\begin{abstract}
		Suppose that a plane wave is incident onto an impenetrable grating  profile of Dirichlet or Impedance type, or a penetrable grating. 
		The grating interface is assumed to be given by a Lipschitz function in two dimensions.	
		We derive stability estimate of the grating diffraction problem via variational method with an explicit dependence of solutions on the incident wavenumber and incident angle. 
		
		\vspace{.2in} {\bf Keywords}: Diffractive gratings, 
		Dirichlet boundary, impedance boundary, transmission conditions, stability estimate.	\end{abstract}
	
	\section{Introduction and main results}
	Diffraction gratings have a long history and are widely used in a wide variety of scientific and technological fields \cite{bao2001,petit1990}. 
	We refer to books \cite{MR4385553,petit1990,MR725334} for  its physical and mathematical background and to \cite{MR1146575,MR1273315,MR1010883,MR1160941} and references therein for mathematical analysis and numerical methods on electromagnetic scattering by diffraction gratings.
	There are two fundamental polarizations in electromagnetic scattering problems. 
	The first is the transverse electric (TE) polarization where the electric field is parallel to the grooves of the grating \cite{MR1342287}.
	As for the second transverse magnetic (TM) polarization, the magnetic field is parallel to the grooves of the grating \cite{MR1417860}. In both TE or TM polarization cases, well-posedness of the grating diffraction problem has been sufficiently studied  under additional conditions imposed on the incident wavenumbers, scattering interfaces and material parameters; see e.g., \cite{arens2010scattering,Kirsch1993,MR1273315,MR1251818,Elschner1998,MR1634859}.
	
	In this paper, we assume that the grating is $ 2\pi$-periodic in the $x_1$-direction and invariant in the 
	$x_3$-direction. In the polarization cases, the three-dimensional scattering problem governed by the Maxwell's equations can be reduced to a two-dimensional diffraction problem modeled by the scalar Helmholtz equation over the $x_1x_2$-plane.
	We investigate	stability estimate of the scattering of time-harmonic (with the time variation $ e^{-i\omega t},\, w>0 $) electromagnetic waves by a perfectly reflective grating or by a penetrable grating in an isotropic lossless medium.	With the help of Rellich identity, we prove well-posedness of the Dirichlet and impedance boundary value problems for impenetrable gratings as well as that  of transmission problems for penetrable gratings. 
	
	
	One feature of our stability estimates is the explicit dependence of the solution on the incident wavenumber $k>0$ and on the incident angle $ \theta$. 
	The wavenumber-dependent estimates were derived in \cite{Chandler2005} for rough surface scattering  due to a compactly supported source term and also in \cite{MR3466840} for acoustic scattering of plane waves by a rectangular cavity. However, to the best of our knowledge, both incident wavenumber- and angle-dependent estimates are not available in the literature.  	
	
	The article is organized as follows. 
	In Section \ref{sec:2} we will give mathematical formulations of the grating diffraction problems and some basic notations. 
	Section \ref{sec:3} is devoted to the variational method for the Dirichlet boundary value problem with a Lipschitz grating profile. We derive a stability result with explicit dependance on the incident angle and also on the wavenumber.	
	In Section \ref{Sec:FM} we discuss the impedance boundary value problem for impenetrable gratings and in Section \ref{sec:TC} the transmission conditions for penetrable gratings, with the same kind of stability estimate. 
	The proof of some preliminary lemmas in Section \ref{sec:3} will be postponed to  Section \ref{sec:App}.
	
	\section{Grating diffraction problems}\label{sec:2}
	Let the profile of the diffraction grating be described by the curve
	$$
	\tilde{\Gamma}:=\{x\in \mathbb{R}^2:x_2=f(x_1),\;x_1 \in \mathbb{R}\},
	$$
	with $f\in C_p^{0,1}$, i.e., $f$ is a $2\pi$-periodic Lipschitz function. 
	We further assume $f(x_1)>f_{-}$ for all $ x_1\in (0,2\pi) $. Denote by $ L>0 $ the Lipschitz constant of $ f $.
	Suppose that the space above $\tilde{\Gamma}$ is filled with a homogeneous and isotropic medium. This implies that the unbounded region 
	$$
	\tilde{\Omega}:=\{x\in \mathbb{R}^2:x_2>f(x_1),\;x_1 \in \mathbb{R}\}
	$$
	is filled with a material whose index of refraction (or wave number) $k$ is a positive constant. Introduce the following notations in one period of the grating profile
	\ben
	\Gamma=\Gamma_f&:=\{x\in \mathbb{R}^2:x_2=f(x_1),\;0<x_1<2\pi\},\\
	\Omega=\Omega_f&:=\{x\in \mathbb{R}^2:x_2>f(x_1),\;0<x_1<2\pi\}.
	\enn
	Introduce the artificial boundary
	$$
	\Gamma_R:=\{(x_1,R):0\leq x_1 \leq 2\pi \},\qquad R>\Gamma_{\max}:=\max_{0\leq t\leq 2\pi}|f(t)|,
	$$
	and the bounded domain
	$$
	\Omega_R=\Omega_{R,f}:=\{x\in \mathbb{R}^2:f(x_1)<x_2<R,\;0<x_1<2\pi\}.
	$$
	\begin{figure}[htp]
		\centering
			\begin{tikzpicture}
		\tikzstyle{medium} = [fill=pink!20,fill opacity=0.8];
		\tikzstyle{grating} = [fill=blue!10,fill opacity=0.8];
		
		\draw[thick,->] (-4.5,0)--(4.5,0)node[right]{$x_{1}$};
		\draw[thick,->] (-4,-1)--(-4,6)node[right]{$x_2$};
		
		\draw[medium](-4,1.5)--(-1.5,1)--(1,2)--(3,1.5)--(3,5)--(-4,5)--(-4,1.5);
		\draw[grating](-4,1.5)--(-1.5,1)--(1,2)--(3,1.5)--(3,0)--(-4, 0)--(-4,1.5);
		\draw[-stealth,blue,thick] (-2.5,5.5) -- (-1.5,4.5) node [above right] {$u^{i}$};
		\foreach \ang in {45,60,75,105,120,135}{
			\draw[-stealth,blue,thick] (\ang:3)--(\ang:3.75);
		}
		\draw[blue] (-2.5,3) node [above] {$u^{s}$};
		\draw[blue] (2,3.5) node [above] {$u^{s}$};
	
		\node at (-4.25,5) {\scriptsize$R$};
		\node at (3,-.25) {\scriptsize$2\pi$};
		\node at (3.5,5) {$\Gamma_{R}$};
		\node at (3.5,1.5) {$\Gamma$};
		\node at (-0.5,2) {$\Omega_{R}$};
	\end{tikzpicture}
		\caption{Geometry of the grating diffraction problem in two dimensions.}
		\label{F}
	\end{figure}
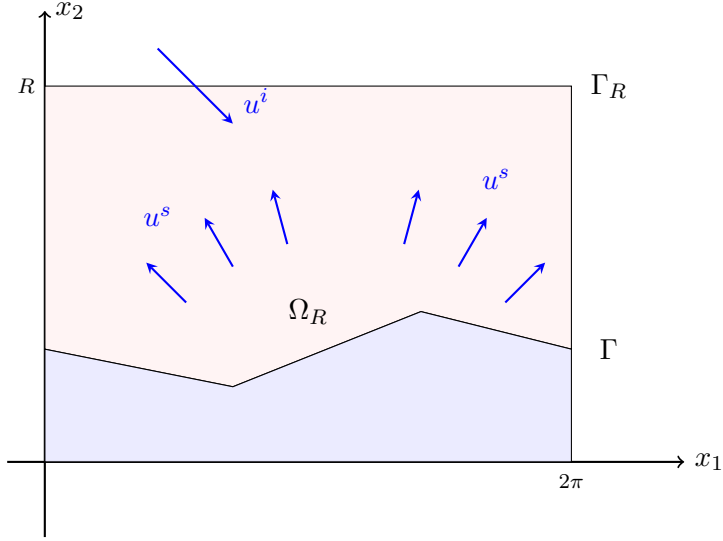 
	We further assume that a plane wave given by
	\be
	u^i(x)=\gamma e^{i\alpha x_1-i\beta x_2},\qquad   \gamma \in \mathbb{C},  \qquad  i=\sqrt{-1}\label{eqn:Ui}
	\en
	is incident onto $\Gamma$ from above, where $\alpha=k\sin{\theta}$, $\beta=k\cos{\theta}$ and $\theta \in \left(-\frac{\pi}{2},\frac{\pi}{2}\right)$ is the incident angle with the positive $ x_ 2 $-axis. 
	Then the scattered field $ u^s=u-u^i $ satisfies the Helmholtz equation 
	\begin{eqnarray}
		\Delta u^s+k^2 u^s=0 \qquad\qquad \rm{in} \quad\Omega. \label{eqn:HelmholtzUs}
	\end{eqnarray}
	Moreover, the scattered field $u^s$ is assumed to be $\alpha$-quasiperiodic in $ x_1 $ in the sense that $ u^s(x_1,x_2)e^{-i\alpha x_1}$ is $2\pi$-periodic in $x_1$. 
	By this definition, we have 
	\begin{eqnarray}
		u^s(x_1+2\pi,x_2)=e^{2i\alpha\pi}u^s(x_1,x_2)\quad \text{ for all } x \in \Omega. \label{eqn:alphaQuasip}
	\end{eqnarray}
	It is obvious that the incident field $u^i$ is $\alpha$-quasiperiodic. 
	Under the assumption (\ref{eqn:alphaQuasip}), the function $ u^s(x_1,x_2)e^{-i\alpha x_1}$ can be expanded into the Fourier series:
	\begin{align*}
		u^s(x_1,x_2)e^{-i\alpha x_1}=\sum_{n \in \mathbb{Z}}u_n(x_2)e^{inx_1},   \qquad  x_ 2>\Gamma_ {\max }.
	\end{align*}
	Substituting $u^s$ into the Helmholtz equation (\ref{eqn:HelmholtzUs}) in $\Omega$ and applying the method of separation of variables, we get the expression of $u^s$ as a sum of plane waves:
	\begin{align*}
		u^s(x_1,x_2)=\sum_{n \in \mathbb{Z}} A_n e^{i\alpha_n x_1+i\beta_n x_2}+B_n e^{i\alpha_n x_1-i\beta_n x_2} ,\qquad A_n,B_n\in \mathbb{C},
	\end{align*}
	where
	\begin{align}
		\alpha_n=n+\alpha,  \qquad
		\begin{array}{l}
			\beta_n=\left\{\begin{array}{ll}
				\sqrt{k^2-|\alpha_n|^2}, & |\alpha_n|\leq k, \\
				i\sqrt{|\alpha_n|^2-k^2}, &  |\alpha_n|> k.
			\end{array}\right. \\ \label{eqn:betaN}
		\end{array} 
	\end{align}
	Physically, the scattered field remains bounded as $x_2 \rightarrow \infty$. 
	Hence, $u^s$ is only composed of bounded outgoing waves in $\Omega$, leading to the well-known Rayleigh expansion condition:
	\begin{eqnarray}
		u^s(x)=\sum_{n\in \mathbb{Z}}u_n e^{i\alpha_n x_1+i\beta_n x_2},\qquad x_2>\Gamma_{\max},  \label{eqn:Us}
	\end{eqnarray}
	with the Rayleigh coefficients $u_n\in \mathbb{C}$.
	Obviously, $u^s$ in (\ref{eqn:Us}) can be split into the finite sum $\sum_{|\alpha_n|\leq k}$ of outgoing plane waves and the infinite sum $\sum_{|\alpha_n|> k}$ of exponentially decaying waves, which are called surface or evanescent waves.
	We refer to Figure \ref{F} for an illustration of the grating diffraction problem. 
	
	We now introduce periodic and quasiperiodic Sobolev spaces to be used in this paper. For $s\in \mathbb{R}$, $s\geq 0$, the Sobolev space $H^s(0,2\pi)$ of periodic functions is defined as the completion of $\left\{u \big|_{[0,2\pi]}:u \mbox{ is trigonometric polynomial}\right\}$ with respect to the inner product
	\begin{align*}
		\langle u,v \rangle:=\sum_{n \in \mathbb{Z}}(k^2+n^2)^s u_n \Bar{v}_n,
	\end{align*}
	where $u_n$ and $v_n$ are Fourier coefficients of $u$ and $v$, respectively. 
	The periodic Sobolev space $H_p^s(\Gamma)$ and the $\alpha$-quasiperiodic Sobolev space $H_{\alpha}^s(\Gamma)$ can be defined, respectively, by
	\begin{align*}
		H_p^s(\Gamma) &:= \left\{u:\Gamma\rightarrow \mathbb{C},\;u(x_1,f(x_1))\in H^s(0,2\pi)\right\},\\
		H_{\alpha}^s(\Gamma) &:= \left\{u:\Gamma\rightarrow \mathbb{C},\;u(x_1,f(x_1))e^{-i\alpha x_1}\in H_p^s(0,2\pi)\right\}.
	\end{align*}
	We also define
	\begin{align*}
		H_{\alpha}^1(\Omega_R)&:=\left\{u\in H^1(\Omega_R):u(x_1,x_2)e^{-i\alpha x_1} \mbox{ is }2\pi\mbox{-periodic with respect to }x_1\right\}.\\
		H_{\rm{loc},\alpha}^1(\Omega)&:=\left\{u:u\in H_{\alpha}^1( \Omega_b)\mbox{ for any } b >\Gamma_{\max}\right \}.
	\end{align*}
	
	\section{Dirichlet boundary value problem}\label{sec:3}	
	In this section, we consider the TE polarization of electromagnetic scattering from a perfectly conducting grating. 
	The problem we wish to analyze is to find $u \in H_{loc,\alpha}^1(\Omega)$ such that
	\begin{align}  
		\begin{array}{l}
			(\mbox{DBVP}):\qquad\left\{\begin{array}{ll}
				\Delta u+k^2u = 0   &\mbox{in} \quad \Omega, \\
				u = 0          &\mbox{on} \quad \Gamma,\\
				u-u^i = \sum_{n \in \mathbb{Z}}u_n e^{i\alpha_n x_1+i\beta_n x_2} & x_2 > \Gamma_{\max}. 
			\end{array}\right. \\ \label{Dirichlet}
		\end{array} 
	\end{align}
	The purpose of this section is to derive a stability estimate of the Dirichlet problem (\ref{Dirichlet}) via the variational approach.
	Define the space
	$$X_R=X(\Omega_{R}):=\{u \in H^1_{\alpha}(\Omega_R):u=0 \; \rm{on} \; \Gamma \},$$ equipped with the norm 
	$$\|u\|_{X_R}^2=k^2\|u\|_{L^2(\Omega_R)}^2+\|\nabla  u\|_{L^2(\Omega_R)}^2.$$
	We define the Dirichlet-to-Neumann map $T:H_{\alpha}^{1/2}(\Gamma_R)\rightarrow H_{\alpha}^{-1/2}(\Gamma_R)$ on the artificial boundary $\Gamma_R$ by
	\begin{align*}
		(Tg)(x_1):=\sum_{n \in \mathbb{Z}}i\beta_n g_n e^{i\alpha_n x_1}, \qquad g(x_1):=\sum_{n \in \mathbb{Z}} g_n e^{i\alpha_n x_1}\in H_{\alpha}^{1/2}(\Gamma_R),
	\end{align*}
	which is equivalent to the Rayleigh expansion (\ref{eqn:Us}).
	The operator $T$ is well defined and bounded because $\beta_n=i|\alpha_n|+O\left(\frac{1}{|n|}\right)$ as $|n|\rightarrow +\infty$. Clearly,
	\begin{align*}
		T\left(u^s \big|_{x_2=R}\right)=\frac{\partial u^s}{\partial x_2} \Big|_{\Gamma_R},\qquad 
		T\left(u^i \big|_{x_2=R}\right)=i\beta e^{-i\beta R}\gamma e^{i\alpha x_1}.
	\end{align*}
	This implies that 
	\begin{align}
		T\left(u^i \big|_{\Gamma_R}\right)-\frac{\partial u^i}{\partial \nu}\Big|_{\Gamma_R}=2i\beta e^{-i\beta R}\gamma e^{i\alpha x_1},\label{lefterm}
	\end{align}
	where $\nu=(0,1)^\top$ is the normal direction on $ \Gamma_ R $ pointing to $\Omega$.
	
	\begin{theorem}\label{th:dependence}
		Let $u^i$ be a plane wave given by (\ref{eqn:Ui}) and suppose that $u \in X_R$ is the unique solution to the Dirichlet problem (\ref{Dirichlet}).
		Suppose that $ f_-+1<\Gamma_{\min}:=\min_{0\leq t\leq 2\pi}|f(t)|$, then
		\begin{align} \label{thD}
			\|u\|_{X_R} \leq 2\sqrt{2\pi} \cos\theta|\gamma |C,
		\end{align} 
		where 
		\begin{align*}
			C&= \sqrt{k M^2+4k^4(R-f_-)^3},\\
			M&= 4k^3(R-f_-)^3+2k^2(R-f_-)^2+4k^2(R-f_-)^3k\cos{\theta}+1.
		\end{align*}
	\end{theorem}
	\begin{remark}
		The stability estimate (\ref{thD}) obviously yields uniqueness of ($ \mbox{DBVP} $).
		If the incident angle $ \theta  $ tends to $ \frac{\pi}{2}  $ or $ -\frac{\pi}{2}  $, the unique solution will converge to zero. 
		At the normal incidence (i.e. $ \theta= 0 $), the energy of the solution is larger than other incident angles. 
	\end{remark}
	It is easy to get the following variational formulation in one periodic cell:
	find $ u\in X_R $ such that
	\begin{equation}
		a(u,v)=F(v)  \quad\text{ for all } v \in X_R   \label{eqn:varitional formula},
	\end{equation}
	where
	
	\begin{align}
		a(u,v)&:=\int_{X_R}\nabla u \cdot \nabla \Bar{v}-k^2u\Bar{v}\, dx-\int_{\Gamma_R}(Tu)\Bar{v}\, ds, \label{eqn:dirichletAuv}\\
		F(v)&:=-2i\beta e^{-i\beta R}\gamma \int_0^{2\pi}e^{i\alpha x_1}\Bar{v}(x_1,R)\, dx_1. \no 
	\end{align}
	Noting that on $ \Gamma_{R} $, we have used (\ref{lefterm}) and the following identity:
	\ben
	\partial_{\nu}u=\partial_{\nu}u^s+\partial_{\nu}u^i=Tu^s+\partial_{\nu}u^i=Tu-Tu^i+\partial_{\nu}u^i.
	\enn
	We also need the following estimates.
	\begin{lemma}
		If $u\in X_R$ is a solution to (\ref{Dirichlet}), then we have
		\begin{eqnarray}
			\left\|\frac{\partial u}{\partial x_2}\right\|_{L^2(\Omega_R)}^2 \leq
			(R-f_-)\int_{\Gamma_R}\left|\frac{\partial u}{\partial x_2}\right|^2-\left|\frac{\partial u}{\partial x_1}\right|^2+k^2|u|^2\,ds+\int_{\Omega_R}|\nabla u|^2-k^2|u|^2\,dx .\label{eqn:estiD2u}
		\end{eqnarray}
	\end{lemma}
	\begin{proof}
		For the grating profile function $ f\in C_p^{0,1} $, there exists a sequence of infinitely smooth  profile functions $ f_j,\,j\in \mathbb{N} $ such that $ \|f_j-f\|_{C_p^{0,1}}\to\infty $ as $ j \to\infty$. 
		Define 
		$$ \Omega_{R,j}:=\{x\in \mathbb{R}^2:f_j(x_1)<x_2<R,\;0<x_1<2\pi\}, $$ 
		and let $ u_j \in X(\Omega_{R,j})\subset X(\Omega_R) $ be the unique solution to the Dirichlet problem (\ref{Dirichlet}). By Remark 3.2 in \cite{Elschner2002}, we have $ u_j \to u,\,j \to\infty$ in $ X_R $. Using the Rellich identity (see Corollary \ref{RemarkrellichNohelmholtz} with $ c=f_- $ in the Appendix) and the fact that $ \nu_2 \geq 0 $ on $ \Gamma $, we obtain
		\begin{align}
			\left\|\frac{\partial u_j}{\partial x_2}\right\|_{L^2(\Omega_{R,j})}^2 
			&\leq  \int_{\Gamma}(x_2-f_-)\nu_2 \left|\frac{\partial u_j}{\partial \nu}\right|^2\,ds+2\int_{\Omega_{R,j}}\left|\frac{\partial u_j}{\partial x_2}\right|^2\,dx \nonumber\\
			&= (R-f_-)\int_{\Gamma_R}\left|\frac{\partial u_j}{\partial x_2}\right|^2-\left|\frac{\partial u_j}{\partial x_1}\right|^2+k^2|u_j|^2\,ds+\int_{\Omega_{R,j}}|\nabla u_j|^2-k^2|u_j|^2\,dx. \label{eqn:estiD2uj}
		\end{align}
		Hence passing the limit $ j \to\infty$ in the inequality (\ref{eqn:estiD2uj}), we obtain (\ref{eqn:estiD2u}). 
		This completes the proof.
	\end{proof}
	The proof of Lemmas \ref{eqn:estimator1} and \ref{eqn:estimator2} below will be given in the Appendix.
	\begin{lemma}\label{eqn:estimator1}
		If $v\in X_R$, then
		\begin{eqnarray*}
			\sqrt{2\pi}\|v\|_{H_{\alpha}^{1/2}(\Gamma_R)}\leq\|v\|_{X_R}.   
		\end{eqnarray*}
	\end{lemma}
	
	\begin{lemma}\label{eqn:estimator2}
		If $v \in H_{\alpha}^1(\Omega_R)$, then
		\begin{eqnarray*}
			\|v\|^2_{L^2(\Omega_R)}\leq(R-f_{-})^2\left\|\frac{\partial v}{\partial x_2}\right\|^2_{L^2(\Omega_R)}+2(R-f_-)\|v\|^2_{L^2(\Gamma)}.   
		\end{eqnarray*}
	\end{lemma}
	From the Rayleigh expansion (\ref{eqn:Us}), we may rewrite the restriction of $ u  $ to $ \Gamma_ R  $ as 
	\begin{eqnarray}
		u|_{\Gamma_R} = \sum_{n \in \mathbb{Z}}u_n e^{i\alpha_n x_1+i\beta_n R}+\gamma e^{i(\alpha x_1-\beta R)}=\sum_{n \in \mathbb{Z}}\tilde{u}_ne^{i\alpha_n x_1}, \label{Qun}
	\end{eqnarray}
	where
	\begin{align*}
		\begin{array}{l}
			\tilde{u}_n:=\left\{\begin{array}{ll}
				u_n e^{i\beta_n R}, & n\neq 0, \\
				u_0 e^{i\beta R}+\gamma  e^{-i\beta R}, & n=0.
			\end{array}\right. \\
		\end{array} 
	\end{align*}
	Using these notations, we have the following lemma.
	\begin{lemma}\label{eqn:right side of rellich}
		Let $u=u^i+u^s\in X_R$ be a solution to (\ref{Dirichlet}). Then we have
		\begin{eqnarray*}
			\int_{\Gamma_R}\left|\frac{\partial u}{\partial x_2}\right|^2-\left|\frac{\partial u}{\partial x_1}\right|^2+k^2|u|^2\, ds \leq 2k{\rm{Im}}\,\int_{\Gamma_R}(Tu)\Bar{u}\, ds-8\pi {\rm{Re}}\,(u_0\bar{\gamma }e^{i2\beta R})|\beta|^2. 
		\end{eqnarray*}
	\end{lemma}
	\begin{proof}
		Recall that
		$$ u=\sum_{n \in \mathbb{Z}}u_n e^{i\alpha_n x_1+i\beta_n x_2}+\gamma e^{i(\alpha x_1-\beta x_2)},\quad x_2>\Gamma_{\max}. $$ 
		Using (\ref{Qun}), we have 
		\begin{align}
			&\quad \int_{\Gamma_R}\left|\frac{\partial u}{\partial x_2}\right|^2-\left|\frac{\partial u}{\partial x_1}\right|^2+k^2|u|^2\, ds \no\\
			&= \int_0^{2\pi} \left|\sum_{n \in \mathbb{Z}} u_n i\beta_n e^{i(\alpha_n x_1+\beta_n R)}-\gamma i\beta e^{i(\alpha x_1-\beta R)} \right|^2 -\left|\sum_{n \in \mathbb{Z}}  i\alpha_n \tilde{u}_ne^{i\alpha_n x_1}\right|^2+k^2\left|\sum_{n \in \mathbb{Z}} \tilde{u}_ne^{i\alpha_n x_1} \right|^2\,dx_1\no\\
			&= \int_0^{2\pi}\left|\sum_{n \neq 0}i \beta_n \tilde{u}_ne^{i\alpha_n x_1}+i\beta e^{i\alpha x_1} (u_0 e^{i\beta R}-\gamma  e^{-i\beta R}) \right|^2\,dx_1  -2\pi\sum_{n \in \mathbb{Z}}|\alpha_n|^2|\tilde{u}_n|^2+2\pi k^2\sum_{n \in \mathbb{Z}}|\tilde{u}_n|^2\no\\
			&\leq 2\pi \left(\sum_{n \neq 0}|\beta_n|^2|\tilde{u}_n|^2+|\beta|^2|u_0 e^{i \beta R}-\gamma  e^{-i\beta R}|^2 \right) -2\pi\sum_{n \in \mathbb{Z}}|\alpha_n|^2|\tilde{u}_n|^2+2\pi k^2\sum_{n \in \mathbb{Z}}|\tilde{u}_n|^2\no\\
			&= 2\pi \sum_{n \in \mathbb{Z}}(|\beta_n|^2-|\alpha_n|^2+k^2)|\tilde{u}_n|^2+ 2\pi|\beta|^2(|u_0 e^{i \beta R}-\gamma  e^{-i\beta R}|^2-|\tilde{u}_0|^2 ).\label{eqn:dirichletDirectesti1}
		\end{align}
		We next calculate the first and second terms on the right-hand side of (\ref{eqn:dirichletDirectesti1}).
		Recalling the definition of $\beta_n$ in (\ref{eqn:betaN}), 
		\ben
		\sum_{n \in \mathbb{Z}}(|\beta_n|^2-|\alpha_n|^2+k^2)|\tilde{u}_n|^2=2\sum_{|\alpha_n| \leq k}|\beta_n|^2|\tilde{u}_n|^2.
		\enn
		By the definition of $\tilde{u}_0$, we have
		\begin{align*}
			|u_0 e^{i \beta R}-\gamma  e^{-i\beta R}|^2-|\tilde{u}_0|^2
			&= |u_0 e^{i \beta R}-\gamma  e^{-i\beta R}|^2-|u_0 e^{i \beta R}+\gamma  e^{-i\beta R}|^2\\
			&= (|u_0 e^{i \beta R}|^2+|\gamma  e^{-i\beta R}|^2-2\mbox{Re}\,(u_0 e^{i \beta R}\bar{\gamma } e^{i \beta R})) \\
			&\quad-(|u_0 e^{i \beta R}|^2+|\gamma  e^{-i\beta R}|^2+2\mbox{Re}\,(u_0 e^{i \beta R}\bar{\gamma } e^{i \beta R}))\\
			&= -4\mbox{Re}\,(u_0 \bar{\gamma } e^{2i \beta R}).
		\end{align*}
		Inserting the previous two equations into (\ref{eqn:dirichletDirectesti1}) gives
		\begin{align}
			\int_{\Gamma_R}\left|\frac{\partial u}{\partial x_2}\right|^2-\left|\frac{\partial u}{\partial x_1}\right|^2+k^2|u|^2\, ds
			&\leq 4\pi \sum_{|\alpha_n| \leq k}|\beta_n|^2|\tilde{u}_n|^2+2\pi |\beta|^2[-4\mbox{Re}\,(u_0 \bar{\gamma }e^{2i\beta R})]\no \\
			&= 4\pi\sum_{|\alpha_n| \leq k}|\beta_n|^2|\tilde{u}_n|^2-8\pi |\beta|^2\mbox{Re}\,(u_0 \bar{\gamma }e^{2i\beta R}).\label{Dleft}
		\end{align}
		Furthermore, by the definition of the Dirichlet-to-Neumann map $T$ and the definition of $ \beta_n $, we have
		\be
		\int_{\Gamma_R}(Tu)\Bar{u}\, ds=\int_{0}^{2\pi}\left(\sum_{n \in \mathbb{Z}}i\beta_n\tilde{u}_ne^{i\alpha_n x_1}\right)\left(\sum_{m \in \mathbb{Z}}\bar{\tilde{u}}_m e^{-i\alpha_m x_1}\right)\,dx_1=2\pi \sum_{n \in \mathbb{Z}}i \beta_n|\tilde{u}_n|^2,\label{Dt}
		\en	
		and therefore
		\be
		\mbox{Im}\,\int_{\Gamma_R}(Tu)\Bar{u}\, ds=2\pi \sum_{|\alpha_n| \leq k}\beta_n|\tilde{u}_n|^2.\label{Dimt}
		\en
		Consequently, combining (\ref{Dleft}) with (\ref{Dimt}) yields
		\begin{align*}
			&\quad 2k\mbox{Im}\,\int_{\Gamma_R}(Tu)\Bar{u}\, ds-8\pi \mbox{Re}\,(u_0\bar{\gamma }e^{i2\beta R})|\beta|^2-\int_{\Gamma_R}\left|\frac{\partial u}{\partial x_2}\right|^2-\left|\frac{\partial u}{\partial x_1}\right|^2+k^2|u|^2\, ds \\
			&\geq 4k\pi \sum_{|\alpha_n| \leq k}\beta_n|\tilde{u}_n|^2-4\pi\sum_{|\alpha_n| \leq k}|\beta_n|^2|\tilde{u}_n|^2 \\
			&= 4\pi \sum_{|\alpha_n| \leq k}\beta_n(k-\beta_n)|\tilde{u}_n|^2\\
			&\geq 0,
		\end{align*}
		which completes the proof.
	\end{proof}
	\textbf{Proof of Theorem \ref{th:dependence}.}
	Choose $v=u$ in the variational formula (\ref{eqn:varitional formula}) and then take the real and imaginary parts respectively to get the following formula
	\begin{eqnarray}
		\int_{\Omega_R}|\nabla u|^2-k^2|u|^2\, dx-\mbox{Re}\,\int_{\Gamma_R}(Tu)\Bar{u}\, ds=\mbox{Re}\,F(u), \label{eqn:real part}\\   
		-\mbox{Im}\,\int_{\Gamma_R}(Tu)\Bar{u}\, ds=\mbox{Im}\,F(u). \label{eqn:image part}
	\end{eqnarray}
	By (\ref{Dt}), we can easily get
	\begin{eqnarray}
		\mbox{Re}\,\int_{\Gamma_R}(Tu)\Bar{u}\, ds=-2\pi\sum_{|\alpha_n|> k}
		|\beta_n||\tilde{u}_n|^2\leq 0,\label{eqn:realLuu} 
	\end{eqnarray}
	Using equation (\ref{eqn:image part}) and Lemma \ref{eqn:right side of rellich}, we can bound the left-hand side of the inequality in Lemma \ref{eqn:right side of rellich} by
	\be
	\int_{\Gamma_R}\left|\frac{\partial u}{\partial x_2}\right|^2-\left|\frac{\partial u}{\partial x_1}\right|^2+k^2|u|^2\, ds \leq -2k\mbox{Im}\,F(u)-8\pi \mbox{Re}\,(u_0\bar{\gamma }e^{i2\beta R})|\beta|^2.  \label{eqn:Newright side of rellich}
	\en
	Inserting equations (\ref{eqn:real part}) and (\ref{eqn:Newright side of rellich}) into (\ref{eqn:estiD2u}) and using (\ref{eqn:realLuu}), we get
	
	\be
	\left\|\frac{\partial u}{\partial x_2}\right\|_{L^2(\Omega_R)}^2
	&\leq& (R-f_-)\int_{\Gamma_R}\left|\frac{\partial u}{\partial x_2}\right|^2-\left|\frac{\partial u}{\partial x_1}\right|^2+k^2|u|^2\,ds+\int_{\Omega_R}|\nabla u|^2-k^2|u|^2\,dx \nonumber\\
	&\leq& (R-f_-)\left[-2k\mbox{Im}\,F(u)-8\pi \mbox{Re}\,(u_0\bar{\gamma }e^{2i\beta R})|\beta|^2\right]\nonumber\\
	&\;&+\mbox{Re}\,\int_{\Gamma_R}(Tu)\Bar{u}\,ds+\mbox{Re}\,F(u) \nonumber\\
	&\leq& (R-f_-)2k|F(u)|+|F(u)|+8\pi(R-f_-)|\beta|^2|\mbox{Re}\,(\bar{\gamma } u_0e^{2i\beta R})|  \nonumber\\
	&\leq& \left[2k(R-f_-)+1\right]|F(u)|+8\pi(R-f_-)|\beta|^2(|\gamma ||\tilde{u}_0|+|\gamma |^2). \label{eqn:dirichletPartialUX2}
	\en
	In the last term of the above inequality, we have used the following identity,
	\begin{align*}
		\bar{\gamma }\tilde{u}_0 e^{i\beta R}=  \bar{\gamma }\left(u_0e^{2i\beta R}+\gamma \right)=\bar{\gamma } u_0e^{2i\beta R}+|\gamma |^2,
	\end{align*}
	which implies,
	\begin{align*}
		\left|\mbox{Re}\,\left(\bar{\gamma } u_0e^{2i\beta R}\right)\right|\leq \left|\mbox{Re}\,\left(\bar{\gamma }\tilde{u}_0 e^{i\beta R}\right)-|\gamma |^2 \right|\leq |\gamma ||\tilde{u}_0|+|\gamma |^2.
	\end{align*}
	On the other hand, since
	\ben
	\Bar{u}\big|_{\Gamma_R}=\sum_{n\in \mathbb Z}\Bar{\tilde{u}}_n e^{-i \alpha_n x_1},\quad \alpha_n=n+\alpha,
	\enn
	we get
	\be
	|F(u)| 
	&=&\left|-2i\beta e^{-i\beta R}\gamma \int_0^{2\pi}e^{i\alpha x_1}\Bar{u}(x_1,R)\,dx_1\right| \nonumber\\
	&=&\left|2\beta \gamma \int_0^{2\pi}\sum_{n \in \mathbb{Z}}\Bar{\tilde{u}}_ne^{i(\alpha-\alpha_n)x_1}\,dx_1\right| \nonumber\\
	&=&|4\pi\beta \gamma \tilde{u}_0|. \label{eqn:DirichletFuVsU0}
	\en
	Hence, (\ref{eqn:dirichletPartialUX2}) can be estimated by
	\be
	\left\|\frac{\partial u}{\partial x_2}\right\|_{L^2(\Omega_R)}^2
	&\leq&\left[2k(R-f_-)+1\right]4\pi|\beta||\gamma ||\tilde{u}_0|+8\pi(R-f_-)|\beta|^2|\gamma ||\tilde{u}_0|\no\\
	&\;& +8\pi(R-f_-)|\beta|^2|\gamma |^2\no\\
	&=& C_1|\tilde{u}_0||\gamma |+8\pi(R-f_-)|\beta|^2|\gamma |^2,\label{eqn:DiriSecondPartialX2}
	\en
	where 
	\begin{eqnarray*}
		C_1=4\pi |\beta|\left[2k(R-f_-)+2|\beta|(R-f_-)+1\right].
	\end{eqnarray*}	
	With the help of (\ref{eqn:real part}) and (\ref{eqn:realLuu}), we get
	\be
	\|u\|_{X_R}^2 
	&=&\|\nabla u\|_{L^2(\Omega_R)}^2+k^2\| u\|_{L^2(\Omega_R)}^2 \nonumber\\
	&=& 2k^2\| u\|_{L^2(\Omega_R)}^2+\mbox{Re}\,\int_{\Gamma_R}(Tu)\Bar{u}\,ds+\mbox{Re}\,F(u) \nonumber\\
	&\leq& 2k^2\| u\|_{L^2(\Omega_R)}^2+|F(u)|.\label{eqn:DirichletUX}
	\en
	Applying Lemma \ref{eqn:estimator2} (note that $ u=0 \;\mbox{on} \; \Gamma$), (\ref{eqn:DirichletFuVsU0}) and (\ref{eqn:DiriSecondPartialX2}), we continue to estimate $\|u\|_{X_R}^2$ in (\ref{eqn:DirichletUX}) by
	\be
	\|u\|_{X_R}^2 
	&\leq& 2k^2(R-f_-)^2\left\|\frac{\partial u}{\partial x_2}\right\|_{L^2(\Omega_R)}^2+4\pi|\beta||\gamma ||\tilde{u}_0| \nonumber\\
	&\leq& 2k^2(R-f_-)^2\left[C_1|\gamma ||\tilde{u}_0|+8\pi(R-f_-)|\beta|^2|\gamma |^2\right]+4\pi|\beta|\gamma ||\tilde{u}_0| \nonumber\\
	&=& C_2|\gamma ||\tilde{u}_0|+C_3|\gamma |^2, \label{eqn:DirichletUX1}
	\en
	where 
	\begin{align*}
		C_3&= 16\pi k^2(R-f_-)^3|\beta|^2,\\
		C_2&= 2k^2(R-f_-)^2C_1+4\pi|\beta|=4\pi |\beta|M,\\
		M&:=4k^3(R-f_-)^3+2k^2(R-f_-)^2+4k^2(R-f_-)^3|\beta|+1.
	\end{align*}	
	Recall Young's inequality
	\be
	C_2|\gamma ||\tilde{u}_0|\leq C_2\left(\epsilon|\tilde{u}_0|^2+\frac{|\gamma |^2}{4\epsilon}\right), \qquad \epsilon >0 \label{eqn:DestimateC2U0C0} 
	\en
	for any $ \epsilon>0 $.
	Using the definition of $\|u\|_{H_{\alpha}^{1/2}(\Gamma_R)}$ and Lemma \ref{eqn:estimator1}, we have
	\be
	\|u\|_{X_R}^2 \geq	2\pi\|u\|_{H_{\alpha}^{1/2}(\Gamma_R)}^2=2\pi\sum_{n\in\mathbb{Z}}(k^2+\alpha_n^2)^{1/2}|\tilde{u}_n|^2\geq 2\pi k|\tilde{u}_0|^2. \label{eqn:U0UX}
	\en	
	Combining (\ref{eqn:DestimateC2U0C0}) and (\ref{eqn:U0UX}) gives
	$$C_2|\gamma ||\tilde{u}_0|\leq \frac{C_2\epsilon}{2\pi k}\|u\|_{X_R}^2+\frac{C_2}{4\epsilon}|\gamma |^2.$$
	Thus, (\ref{eqn:DirichletUX1}) becomes
	\begin{equation*}
		\|u\|_{X_R}^2\leq \frac{C_2\epsilon}{2\pi k}\|u\|_{X_R}^2+\frac{C_2}{4\epsilon}|\gamma |^2+C_3|\gamma |^2,
	\end{equation*}
	which implies
	\begin{eqnarray*}
		\|u\|_{X_R}^2\leq \frac{2\pi k\left(\frac{C_2}{4\epsilon}+C_3\right)}{2\pi k-C_2\epsilon}|\gamma |^2.
	\end{eqnarray*}
	Taking $\epsilon=\frac{\pi k}{C_2}$, we get 
	\begin{eqnarray*}
		\|u\|_{X_R}^2 \leq \left(\frac{C_2}{2\epsilon}+2C_3\right)|\gamma |^2
		=\left(\frac{C_2^2}{2\pi k}+2C_3\right)|\gamma |^2.
	\end{eqnarray*}
	Furthermore, by $\beta=k \cos{\theta}$ and the definitions of $C_2$ and $C_3$, we have
	\begin{eqnarray*}
		\frac{C_2^2}{2\pi k}+2C_3
		= 8\pi k\cos^2\theta M^2+32\pi k^4(R-f_-)^3\cos^2\theta
		= 8\pi \cos^2\theta C^2,
	\end{eqnarray*}
	where $ \epsilon> 0 $ is given as Theorem  \ref{th:dependence}.
	$\hfill\Box$
	
	\section{Impedance boundary value problem}\label{Sec:FM}
	In this section, we suppose that the grating surface is coated by a thin layer of material with the positive surface impedance $\lambda$. 
	Consider the impedance boundary value problem
	\begin{align}  
		\begin{array}{l}
			(\mbox{IBVP}):\qquad\left\{\begin{array}{ll}
				\Delta u+k^2u = 0    &\mbox{in} \quad \Omega, \\
				\partial_{\nu} u+i\lambda u= 0          &\mbox{on} \quad \Gamma,\\
				u-u^i = \sum_{n \in \mathbb{Z}}u_n e^{i\alpha_n x_1+i\beta_n x_2} & x_2 > \Gamma_{\max}. 
			\end{array}\right. \\ \label{Impedance}
		\end{array} 
	\end{align}
	Here $\nu$ is the normal direction at $ \Gamma $ pointing to $\Omega$. We keep the notations used in Section \ref{sec:3}, except for the new definition $ X_R:=H^1_{\alpha}(\Omega_R) $.
	From \cite{Kirsch1993} we know that the problem (\ref{Impedance}) is uniquely solvable in $H_{\alpha}^1(\Omega_R)$. 
	The main stability estimate for a plane wave incidence is stated below.
	\begin{theorem}\label{th:Impedancedependence}
		Let $u^i$ be  a plane wave given by (\ref{eqn:Ui}) and suppose that $u \in H^1_{\alpha}(\Omega_R)$ is the unique solution to (\ref{Impedance}). 
		Suppose that $ R- 1>\Gamma_{\max} $ and $ f_-+1<\Gamma_{\min} $, then 
		\begin{align*}
			\|u\|_{H^1_{\alpha}(\Omega_R)} \leq  2\sqrt{2\pi}  \cos{\theta}|\gamma |C^{*} , 
		\end{align*}
		where
		\begin{align*}
			C^*&=\sqrt{ k \left( 1+\frac{4k^2\tilde{C}^2}{\lambda}\right)^2+8k^4\tilde{C}^2},\\
			\tilde{C}^2&=\frac{4(R-f_-)^2+(2k+1)(R-f_-)^3[2k(R-f_-)+1]}{2\min\left \{\frac{R-f_--1}{(2k+1)(R-f_-)^3},\frac{1}{\sqrt{1+L^2}}\right\}}.
		\end{align*}
		Here $L$ is the Lipschitz constant of $f$.
	\end{theorem}
	Before proving the theorem, we first give the variational formulation of the problem. 
	It is the same as (\ref{eqn:varitional formula}), except that the right-hand side of (\ref{eqn:dirichletAuv}) is subtracted by $\int_{\Gamma}i\lambda u\Bar{v}\,ds$.
	That is
	\begin{equation*}
		a(u,v)=F(v) \quad \text{ for all }  v \in H^1_{\alpha}(\Omega_R)   
	\end{equation*}
	where
	\begin{align*}
		a(u,v)&=\int_{\Omega_R}\nabla u \cdot \nabla \Bar{v}-k^2u\Bar{v}\, dx-\int_{\Gamma}i\lambda u\Bar{v}\,ds-\int_{\Gamma_R}(Tu)\Bar{v}\, ds, \\
		F(v)&=-2i\beta e^{-i\beta R}\gamma \int_0^{2\pi}e^{i\alpha x_1}\Bar{v}(x_1,R)\, dx_1. 
	\end{align*}
	Again we take the real and imaginary parts of the variational formula with $ v=u $ to get 
	\begin{eqnarray}
		\int_{\Omega_R}|\nabla u|^2-k^2|u|^2\, dx-\mbox{Re}\,\int_{\Gamma_R}(Tu)\Bar{u}\, ds=\mbox{Re}\,F(u),    \label{eqn:realUImpedance}\\
		-\int_{\Gamma}\lambda |u|^2\,ds-\mbox{Im}\,\int_{\Gamma_R}(Tu)\Bar{u}\, ds=\mbox{Im}\,F(u). \label{eqn:imageUImpedance}
	\end{eqnarray}
	Analogously to the Dirichlet problem, it holds that 
	\begin{align}
		\mbox{Re}\,\int_{\Gamma_R}(Tu)\Bar{u}\, ds&=-2\pi\sum_{|\alpha_n|> k}
		|\beta_n||\tilde{u}_n|^2 \leq 0, \label{eqn:realLuuImpedance}\\
		\mbox{Im}\,\int_{\Gamma_R}(Tu)\Bar{u}\, ds&=2\pi\sum_{|\alpha_n|\leq k}
		\beta_n|\tilde{u}_n|^2 \geq 0.  \label{eqn:imageLuuImpedance}
	\end{align}
	Using (\ref{eqn:realUImpedance}) and (\ref{eqn:realLuuImpedance})
	
	\begin{align}
		\int_{\Omega_R}|\nabla u|^2\,dx&=k^2\int_{\Omega_R}|u|^2\,dx+\mbox{Re}\,\int_{\Gamma_R}(Tu)\Bar{u}\,ds+\mbox{Re}\,F(u),\nonumber\\
		&\leq k^2\int_{\Omega_R}|u|^2\,dx+|F(u)| . \label{eqn:ImpedanceNabla}
	\end{align}
	Using (\ref{eqn:imageUImpedance}) and (\ref{eqn:imageLuuImpedance}), we get
	\begin{align*}
		\lambda\int_{\Gamma}|u|^2\,ds &= -\mbox{Im}\,\int_{\Gamma_R}(Tu)\Bar{u}\,ds-\mbox{Im}\,F(u) \\
		&= -2\pi \sum_{|\alpha_n|\leq k}\beta_n|\tilde{u}_n|^2-\mbox{Im}\,F(u) \\
		&\leq |F(u)|,
	\end{align*}
	that is
	\begin{eqnarray}
		\int_{\Gamma}|u|^2\,ds \leq \frac{1}{\lambda} |F(u)|. \label{eqn:estimatUgamma}
	\end{eqnarray}	
	To estimate $\|u\|_{L^2(\Omega_R)}$, we need to consider the following auxiliary problem: find $ w\in H_{-\alpha}^1 (\Omega_R)$ such that
	\begin{align}  
		\begin{array}{l}
			(\mbox{AP}):\qquad\left\{\begin{array}{ll}
				\Delta w+k^2w = \bar{u}   &\mbox{in} \quad \Omega, \\
				w = 0          &\mbox{on} \quad \Gamma,\\
				\hat{T}w = \partial_2 w &\mbox{on} \quad\Gamma_{R},\\
			\end{array}\right. \\ \label{AugW}
		\end{array} 
	\end{align}
	where operator $\hat{T}:H_{-\alpha}^{1/2}(\Gamma_R)\rightarrow H_{-\alpha}^{-1/2}(\Gamma_R)$ is defined as
	\ben
	(\hat{T}g)(x_1):=\sum_{n \in \mathbb{Z}}i\hat{\beta}_n \hat{g}_n e^{i\hat{\alpha}_n x_1}, \qquad g(x_1)=\sum_{n \in \mathbb{Z}} \hat{g}_n e^{i\hat{\alpha}_n x_1}\in H_{-\alpha}^{1/2}(\Gamma_R)
	\enn
	with 
	\begin{eqnarray*}
		\hat{\alpha}_n=-\alpha+n,\quad \begin{array}{l}
			\hat{\beta}_n=\left\{\begin{array}{ll}
				\sqrt{k^2-|\hat{\alpha}_n|^2}, & |\hat{\alpha}_n|\leq k, \\
				i\sqrt{|\hat{\alpha}_n|^2-k^2}, & |\hat{\alpha}_n|> k.
			\end{array}\right. \\
		\end{array}
	\end{eqnarray*}
	\begin{lemma}\label{LemmaW}
		Let $w\in H^1_{-\alpha}(\Omega_R)$ be a solution to (\ref{AugW}). 
		Assume $ f_-+1<\Gamma_{\min} $, then we have
		\begin{eqnarray*}
			\|w\|_{L^2(\Omega_R)}+\|\partial_{\nu}w\|_{L^2(\Gamma)}\leq \tilde{C} \|u\|_{L^2(\Omega_R)},
		\end{eqnarray*}
		where 
		$$\tilde{C}^2
		=\frac{4(R-f_-)^2+(2k+1)(R-f_-)^3[2k(R-f_-)+1]}{2\min\left \{\frac{R-f_--1}{(2k+1)(R-f_-)^3},\frac{1}{\sqrt{1+L^2}}\right\}}.$$
	\end{lemma}
	\begin{proof}	
		First, we assume that $ \Gamma $ is an infinitely smooth curve, so that $ w\in H_{-\alpha}^2(\Omega_{R}) $ and $ \max_{x_1\in (0,2\pi)}|f'(x_1)|<L $. 
		By Rellich identity (see Corollary \ref{RemarkrellichNohelmholtz} with $ c=f_- $ in the Appendix), we get
		
		\begin{eqnarray}
			2\int_{\Omega_R}|\partial_2 w|^2\,dx+\int_{\Gamma}(x_2-f_-)\nu_2|\partial_{\nu}w|^2\,ds 
			=I_1+I_2+I_3,  \label{eqn:rellechW}
		\end{eqnarray}
		where
		\begin{align*}
			I_1&=-2\mbox{Re}\,\int_{\Omega_R}(x_2-f_-)\partial_2\Bar{w} \Bar{u}\,dx, \\
			I_2&=(R-f_-)\int_{\Gamma_R}|\partial_2 w|^2-|\partial_1 w|^2+k^2|w|^2\,ds,\\
			I_3&=\int_{\Omega_R}|\nabla w|^2-k^2|w|^2\,dx.
		\end{align*}
		Below we shall estimate $I_1$, $I_2$ and $I_3$ in (\ref{eqn:rellechW}) separately. 
		First, $ I_1 $ can be bounded by
		\begin{align}
			|I_1|\leq 2(R-f_-)\|\partial_2 w\|_{L^2(\Omega_R)}\|u\|_{L^2(\Omega_R)}. \label{eqn:estimatI1}
		\end{align}
		Expanding $ w $ into the Rayleigh series
		$$w=\sum_{n \in \mathbb{Z}}w_n e^{i(\hat{\alpha}_n x_1+\hat{\beta}_n x_2)},\qquad x_2>\Gamma_{\max},$$
		we find	
		\begin{eqnarray}
			\hat{T}\left(w|_{x_2=R}\right) &=&\sum_{n \in \mathbb{Z}}i \hat{\beta}_n w_n e^{i(\hat{\alpha}_n x_1+\hat{\beta}_n R)}, \no\\
			\int_{\Gamma_R}(\hat{T}w)\Bar{w}\,ds &=&2\pi\sum_{n\in \mathbb{Z}}i \hat{\beta}_n|w_n|^2,\no\\
			\mbox{Re}\,\int_{\Gamma_R}(\hat{T}w)\Bar{w}\,ds &=& -2\pi\sum_{|\hat{\alpha}_n|> k}| \hat{\beta}_n||w_n|^2\leq 0,\label{ReW}\\
			\mbox{Im}\,\int_{\Gamma_R}(\hat{T}w)\Bar{w}\,ds &=&2\pi\sum_{|\hat{\alpha}_n|\leq k} \hat{\beta}_n|w_n|^2 \geq 0,\no\\
			\int_{\Gamma_R}|\partial_2 w|^2-|\partial_1 w|^2+k^2|w|^2\,ds &\leq& 2k \mbox{Im}\,\int_{\Gamma_R}(\hat{T}w)\Bar{w}\,ds. \label{eqn:leftRellichW}
		\end{eqnarray}
		From the variational formulation of (\ref{AugW}), we get
		\begin{eqnarray}
			\int_{\Omega_R}|\nabla w|^2-k^2|w|^2 \,dx-\int_{\Gamma_R}(\hat{T}w)\Bar{w}\,ds=-\int_{\Omega_R}\Bar{u}\Bar{w}\,dx. \label{vfw}
		\end{eqnarray}
		Taking the real and imaginary parts of (\ref{vfw}) gives 
		\begin{align}
			\int_{\Omega_R}|\nabla w|^2-k^2|w|^2 \,dx-\mbox{Re}\,\int_{\Gamma_R}(\hat{T}w)\Bar{w}\,ds&=-\mbox{Re}\,\int_{\Omega_R}\Bar{u}\Bar{w}\,dx, \label{eqn:realW}\\
			\mbox{Im}\,\int_{\Gamma_R}(\hat{T}w)\Bar{w}\,ds&=\mbox{Im}\,\int_{\Omega_R}\Bar{u}\Bar{w}\,dx. \label{eqn:imageW}
		\end{align}
		Combining (\ref{eqn:leftRellichW}) and (\ref{eqn:imageW})
		\begin{eqnarray}
			\int_{\Gamma_R}|\partial_2 w|^2-|\partial_1 w|^2+k^2|w|^2\,ds \leq 2k\mbox{Im}\,\int_{\Omega_R}\Bar{u}\Bar{w}\,dx \leq 
			2k \|u\|_{L^2(\Omega_R)}\|w\|_{L^2(\Omega_R)}.\label{eqn:ImpeRightI2}
		\end{eqnarray}
		Inserting (\ref{eqn:ImpeRightI2}) into $I_2$ gives the estimate of $I_2$:
		\begin{eqnarray}
			|I_2|&\leq& 2k(R-f_-)\|u\|_{L^2(\Omega_R)}\|w\|_{L^2(\Omega_R)}.\label{eqn:estimateI2}
		\end{eqnarray}
		Using (\ref{eqn:realW}) and (\ref{ReW}), we have
		\begin{eqnarray}
			I_3 = \mbox{Re}\,\int_{\Gamma_R}(\hat{T}w)\Bar{w}\,ds-\mbox{Re}\,\int_{\Omega_R}\Bar{u}\Bar{w}\,dx
			\leq \|u\|_{L^2(\Omega_R)}\|w\|_{L^2(\Omega_R)}.\label{eqn:estimateI3}
		\end{eqnarray}	
		Next, we derive a low bound of the left-hand side of (\ref{eqn:rellechW}).
		The unit normal vector to $ \Gamma $ is given by $\nu=(\nu_1,\nu_2)^\top=\frac{(-f'(x_1),1)^\top}{\sqrt{1+|f'(x_1) |^2}}$, implying $\nu_2 \geq C_L$, where $C_L=\frac{1}{\sqrt{1+L^2}}\leq 1$. 
		Owing to $ f( x_ 1)-f_ -\geq 1 $, we have 
		\begin{eqnarray}
			2\int_{\Omega_R}|\partial_2 w|^2\,dx+\int_{\Gamma}(x_2-f_-)\nu_2|\partial_{\nu}w|^2\,ds \geq 2\|\partial_2 w \|_{L^2(\Omega_R)}^2+C_L\|\partial_{\nu} w\|_{L^2(\Gamma)}^2.\label{eqn:ImpedanceLeftRellich}
		\end{eqnarray}
		Substituting (\ref{eqn:ImpedanceLeftRellich}), (\ref{eqn:estimatI1}), (\ref{eqn:estimateI2}), (\ref{eqn:estimateI3}) into (\ref{eqn:rellechW}), we get 
		\begin{align}
			&\quad 2\|\partial_2 w \|_{L^2(\Omega_R)}^2+C_L\|\partial_{\nu} w\|_{L^2(\Gamma)}^2 \nonumber\\
			&\leq 2(R-f_-)\left(\epsilon_1 \|\partial_2 w\|_{L^2(\Omega_R)}^2+\frac{1}{4\epsilon_1} \|u\|_{L^2(\Omega_R)}^2\right)
			+2k(R-f_-)\left( \epsilon_2 \| w\|_{L^2(\Omega_R)}^2\right.\nonumber\\
			&\quad \left.+\frac{1}{4\epsilon_2} \|u\|_{L^2(\Omega_R)}^2 \right)
			+\left( \epsilon_2 \| w\|_{L^2(\Omega_R)}^2+\frac{1}{4\epsilon_2} \|u\|_{L^2(\Omega_R)}^2\right)\nonumber\\
			&\leq 2\epsilon_1(R-f_-)\|\partial_2 w \|_{L^2(\Omega_R)}^2+\epsilon_2[2k(R-f_-)+1]\| w\|_{L^2(\Omega_R)}^2 \nonumber\\
			&\quad +\left(\frac{R-f_-}{2\epsilon_1}+\frac{2k(R-f_-)+1}{4\epsilon_2}\right)\|u\|_{L^2(\Omega_R)}^2 \nonumber\\
			&= C_4(\epsilon_1)\|\partial_2 w \|_{L^2(\Omega_R)}^2+C_5(\epsilon_2)\| w\|_{L^2(\Omega_R)}^2+C_6(\epsilon_1,\epsilon_2)\|u\|_{L^2(\Omega_R)}^2,\label{eqn:estimateRellich1}
		\end{align}
		where
		\begin{eqnarray*}
			C_4 (\epsilon_1)&=& 2\epsilon_1(R-f_-), \\
			C_5(\epsilon_2) &=& \epsilon_2[2k(R-f_-)+1],\\
			C_6(\epsilon_1,\epsilon_2) &=& \frac{R-f_-}{2\epsilon_1}+\frac{2k(R-f_-)+1}{4\epsilon_2},
		\end{eqnarray*}
		and $ \epsilon_1,\epsilon_2>0 $ are arbitrary.
		Rewriting (\ref{eqn:estimateRellich1}), we get
		\be
		(2-C_4(\epsilon_1))\|\partial_2 w\|_{L^2(\Omega_R)}^2-C_5(\epsilon_2)\|w\|_{L^2(\Omega_R)}^2+C_L\|\partial_{\nu} w\|_{L^2(\Gamma)}^2 \leq C_6(\epsilon_1,\epsilon_2)\|u\|_{L^2(\Omega_R)}^2. \label{eqn:ImpEstiRellichNew}
		\en
		Using Lemma \ref{eqn:estimator2} and noting again $ w=0 $ on $ \Gamma $, we have
		\be
		\|\partial_2 w\|_{L^2(\Omega_R)}^2 \geq \frac{1}{(R-f_-)^2}\|w\|_{L^2(\Omega_R)}^2.\label{eqn:ImpeNormW2VsW}
		\en
		Substituting (\ref{eqn:ImpeNormW2VsW}) into (\ref{eqn:ImpEstiRellichNew}), we get
		\begin{eqnarray}
			C_7(\epsilon_1,\epsilon_2) \|w\|_{L^2(\Omega_R)}^2+C_L\|\partial_{\nu} w\|_{L^2(\Gamma)}^2 \leq C_6(\epsilon_1,\epsilon_2)\|u\|_{L^2(\Omega_R)}^2,\label{eqn:estimateRellich2}
		\end{eqnarray}
		where 
		\begin{align*}
			C_7(\epsilon_1,\epsilon_2)=\frac{2-C_4(\epsilon_1)}{(R-f_-)^2}-C_5(\epsilon_2).
		\end{align*}
		Taking appropriate $\epsilon_1$ and $\epsilon_2$ to satisfy $C_7(\epsilon_1,\epsilon_2)>0$ and $C_8(\epsilon_1,\epsilon_2)=\min\{C_7(\epsilon_1,\epsilon_2),C_L\}$, we have from (\ref{eqn:estimateRellich2}) that 
		\begin{align*}
			\frac{C_8(\epsilon_1,\epsilon_2)}{2}\left(\|w\|_{L^2(\Omega_R)}+\|\partial_{\nu} w\|_{L^2(\Gamma)}\right)^2 &\leq C_8(\epsilon_1,\epsilon_2)\left(\|w\|_{L^2(\Omega_R)}^2+\|\partial_{\nu} w\|_{L^2(\Gamma)}^2 \right)\\
			&\leq C_6(\epsilon_1,\epsilon_2)\|u\|_{L^2(\Omega_R)}^2.
		\end{align*}
		That is,
		\begin{eqnarray*}
			\|w\|_{L^2(\Omega_R)}+\|\partial_{\nu} w\|_{L^2(\Gamma)}\leq \tilde{C}(\epsilon_1,\epsilon_2)\|u\|_{L^2(\Omega_R)},
		\end{eqnarray*}
		where $ \tilde{C}(\epsilon_1,\epsilon_2)= \sqrt{\frac{2C_6(\epsilon_1,\epsilon_2)}{C_8(\epsilon_1,\epsilon_2)}} $.
		Next, we analyze the sufficiently small $\epsilon_1$ and $\epsilon_2$ satisfying the constraint $C_7(\epsilon_1,\epsilon_2)>0$. 
		Then we have 
		\begin{align*}
			\frac{2-2\epsilon_1(R-f_-) }{(R-f_-)^2}-\epsilon_2[2k(R-f_-)+1]>0,
		\end{align*}
		which implies that
		\begin{align*}
			\epsilon_2<\frac{2[1-\epsilon_1(R-f_-)]}{(R-f_-)^2[2k(R-f_-)+1]}.
		\end{align*}
		By taking $\epsilon_1=\frac{1}{2(R-f_-)}$, we have
		\begin{align*}
			\epsilon_2<\frac{1}{2k(R-f_-)^3+(R-f_-)^2}.
		\end{align*}
		Recalling that $ R-f_->1 $ and taking $\epsilon_2=\frac{1}{(2k+1)(R-f_-)^3}$, we have
		\begin{align*}
			C_4 &= 1,\\
			C_5 &= \frac{2k(R-f_-)+1}{(2k+1)(R-f_-)^3}, \\
			C_6 &= (R-f_-)^2+\frac{(2k+1)(R-f_-)^3[2k(R-f_-)+1]}{4},\\
			C_7 &= \frac{R-f_--1}{(2k+1)(R-f_-)^3} ,\\
			C_8 &= \min\left \{\frac{R-f_--1}{(2k+1)(R-f_-)^3},\frac{1}{\sqrt{1+L^2}}\right\}.
		\end{align*}
		Next, we start to calculate $\tilde{C} $.
		\begin{align}
			\tilde{C}^2:=	\frac{2C_6}{C_8}=\frac{4(R-f_-)^2+(2k+1)(R-f_-)^3[2k(R-f_-)+1]}{2\min\left \{\frac{R-f_--1}{(2k+1)(R-f_-)^3},\frac{1}{\sqrt{1+L^2}}\right\}}.\label{qC}
		\end{align}
		
		Now we consider the case of a Lipchitz profile function $ f\in C_p^{0,1} $.
		As we know (\ref{AugW}) has a unique solution $ w\in H^1_{-\alpha}(\Omega_R)$, we follow Theorem 3.3 in \cite{Elschner2002}. 
		We choose $ C^{\infty} $ profiles $ \Gamma_{j}=\Gamma_{f_j} $ such that $ \Omega_{R,j}\subset \Omega_R $ for all $ \|f_j-f\|_{C_p^{0,1}}\to 0 $, 
		$ f_- \leq f_j(x_1)\leq R,\text{ for all } j,\; x_1$ and the Lipschitz constant is no more than $ L $. 
		We also assume $ \min_{x_1\in (0,2\pi)}f_j(x_1)-f_- \geq 1 $. 
		Let $ w_j\in H^1_{-\alpha}(\Omega_{R,j})\subset H^1_{-\alpha}(\Omega_{R}) $ be the solution of problem (\ref{AugW}) for $ \Omega_{R,j} $. 
		Applying Theorem 3.1 in \cite{Elschner2002} and the arguments used in the case of an infinitely smooth profile, we obtain
		\be
		\left(\int_{\Omega_{R,j}} |w_j|^2\right)^{1/2}+\left(\int_{\Gamma_{j}}|\partial_\nu w_j|^2\right) ^{1/2}\leq \tilde{C}\left(\int_{\Omega_{R,j}} |u|^2\right)^{1/2}, \label{eqn:estiWj}
		\en
		where $ \tilde{C} $ is given by (\ref{qC}).
		From (\ref{eqn:estiWj}) we then get $ \partial_\nu w_j \rightharpoonup v  $ in $ L^2(0,2\pi) $ (for some subsequence) and by Remark 3.2 in \cite{Elschner2002}, $ w_j \to w $ in $ H^1_{-\alpha}(\Omega_R) $. 
		It remains to check that $ v $ conincides with the $ H_{-\alpha}^{-1/2} $ trace of $ \partial_\nu w $ on $ \Gamma $. 
		Remark 2.4.5 in \cite{Necas1967} implies that $ \phi \big|_{\Gamma_{j}}\to \phi \big|_{\Gamma} $ in $ L^2(0,2\pi) $ for any $ \phi \in H_{-\alpha}^1(\Omega) $. 
		From the variational formulations of $  w_j $ and for all $ \phi \in H_{-\alpha}^1(\Omega) $, we have
		\be
		\int_{\Omega_{R,j}} \nabla w_j\cdot \nabla \bar{\phi}-k^2 w_j \bar{\phi}\,dx-\int_{\Gamma_R}(\hat{T}w_j)\bar{\phi}\,ds+\int_{\Gamma_{j}}(\partial_\nu w_j)\bar{\phi}\,ds=-\int_{\Omega_{R,j}}\bar{u}\bar{\phi}\,dx . \label{lim1}
		\en
		By the variational formulation of (\ref{AugW}), we have
		\be
		\int_{\Omega_{R}} \nabla w\cdot \nabla \bar{\phi}-k^2 w \bar{\phi}\,dx-\int_{\Gamma_R}(\hat{T}w)\bar{\phi}\,ds+\int_{\Gamma}(\partial_\nu w)\bar{\phi}\,ds=-\int_{\Omega_{R}}\bar{u}\bar{\phi}\,dx .\label{lim2}
		\en
		Passing to the limit $ j\to\infty $ of (\ref{lim1}) and using $ \partial_\nu w_j \rightharpoonup v $, together with (\ref{lim2}), we have $ v=\partial_{\nu} w $ on $ \Gamma $.
		Finally, we take $ j\to\infty $ in (\ref{eqn:estiWj}), which completes the proof.
	\end{proof}
	The following result will be used in the proof of Theorem \ref{th:Impedancedependence}.
	\begin{lemma}\label{eqn:ImpedanceWvsW0}
		Assume $ R-1>\Gamma_{\max} $, then for $ w\in H_{-\alpha}^1(\Omega_{R}) $, we have 
		$ |w_0|\leq \frac{1}{ \sqrt{2\pi}}\|w\|_{L^2(\Omega_R)}, $
		where $ w_0=\frac{1}{2\pi}\int_{0}^{2\pi}w(x_1,R)e^{-i\alpha x_1}e^{-i\beta R}\,dx$.
	\end{lemma}
	\begin{proof}
		Define $ D:=(0,2\pi)\times (R-1,R) $, then $ D\subset \Omega_{R} $ and
		\begin{align*}
			\int_D|w|^2\,dx&=\int_0^{2\pi}\int_{R-1}^R|w|^2\,dx_2dx_1\\
			&= \int_0^{2\pi}\int_{R-1}^R \sum_{n\in \mathbb{Z}}w_n e^{i(\hat{\alpha}_n x_1+\hat{\beta}_n x_2)}\sum_{m \in \mathbb{Z}} \bar{w}_m e^{-i\hat{\alpha}_m x_1-i\bar{\hat{\beta}}_m x_2} \,dx_2dx_1\\
			&=2\pi\sum_{n\in \mathbb{Z}}|w_n|^2\int_{R-1}^R e^{i(\hat{\beta}_n -\bar{\hat{\beta}}_n)x_2}\,dx_2 \\
			&=2\pi\left [\sum_{|\hat{\alpha}_n|\leq k}|w_n|^2+\sum_{|\hat{\alpha}_n|> k}|w_n|^2\frac{1}{2|\hat{\beta}_n|}(e^{-2|\hat{\beta}_n|(R-1)}-e^{-2|\hat{\beta}_n|R}) \right] \geq 2\pi |w_0|^2,\\
		\end{align*}
		Therefore
		\begin{eqnarray*}
			\|w\|_{L^2(\Omega_R)} \geq \|w\|_{L^2(D)} \geq \sqrt{2\pi}|w_0|,
		\end{eqnarray*}
		which completes the proof.
	\end{proof}
	\textbf{Proof of Theorem \ref{th:Impedancedependence}.}
	Let $ w \in H^1_{\alpha}(\Omega_{R}) $ be the unique solution to (\ref{AugW}).
	Using Green's formula, the Helmholtz equation for $u=u^i+u^s$ and the fact that $w=0\;\rm{on} \;\Gamma$, we have
	\begin{align}
		\|u\|_{L^2(\Omega_R)}^2 
		&=\int_{\Omega_R}u\Bar{u}\,dx  \nonumber\\
		&=\int_{\Omega_R}u(\Delta w+k^2 w)\,dx \nonumber\\
		&=\int_{\partial\Omega_R}\left(\partial_{\nu}wu-\partial_{\nu}u w\right)\,ds+\int_{\Omega_R}w(\Delta u+k^2 u)\,dx  \nonumber\\
		&= \int_{\Gamma_R}\left(\partial_{\nu}wu-\partial_{\nu}u w\right)\,ds+\int_{\Gamma}\left(\partial_{\nu}wu-\partial_{\nu}uw\right)\,ds  \nonumber\\
		&= \int_{\Gamma_R}\left[\partial_{\nu}w(u^i+u^s)-\partial_{\nu}(u^i+u^s) w\right]\,ds-\int_{\Gamma}\partial_{\nu}wu\,ds \nonumber\\
		&= \int_{\Gamma_R}\left(\partial_{\nu}wu^s-\partial_{\nu}u^s w\right) \,ds+\int_{\Gamma_R}\left(\partial_{\nu}wu^i-\partial_{\nu}u^i w\right) \,ds-\int_{\Gamma}\partial_{\nu}w u\,ds.\label{eqn:impedanceEstimateU}
	\end{align}
	The first integral on the right-hand side of (\ref{eqn:impedanceEstimateU}) vanishes, because
	\begin{align}
		\int_{\Gamma_R}\left(\partial_{\nu}w u^s-\partial_{\nu}u^s w\right) \,ds 
		&=  \int_{0}^{2\pi}\left(\sum_{n\in \mathbb{Z}}i \hat{\beta}_n w_n e^{i(\hat{\alpha}_n x_1+\hat{\beta}_n R)} \sum_{m\in \mathbb{Z}}  u_m e^{i(\alpha_m x_1+\beta_m R)} \right. \nonumber\\
		&\quad \left.-\sum_{n\in \mathbb{Z}}i \beta_n u_n e^{i(\alpha_n x_1+\beta_n R)}\sum_{m\in \mathbb{Z}} w_m e^{i(\hat{\alpha}_m x_1+\hat{\beta}_m R)}\right)\,dx_1 \nonumber\\
		&=  \int_{0}^{2\pi}\left(\sum_{m,n\in \mathbb{Z}}i \hat{\beta}_n w_n  u_m e^{i[(-\alpha+n)+(\alpha+m)] x_1+i(\hat{\beta}_n+\beta_m) R} \right. \nonumber\\
		&\quad \left.-\sum_{n\in \mathbb{Z}}i \beta_n u_n w_m e^{i[(\alpha+n)+(-\alpha+m)] x_1+i(\beta_n+\hat{\beta}_m) R}\right)\,dx_1 \nonumber\\
		&= 2\pi\sum_{n \in \mathbb{Z}}i\hat{\beta}_n w_n u_{-n}e^{i2\hat{\beta}_nR}-2\pi\sum_{n \in \mathbb{Z}}i\beta_n u_n w_{-n}e^{i2\beta_n R}\nonumber\\
		&= 2\pi\sum_{n \in \mathbb{Z}}\left(i\beta_n u_n w_{-n}e^{i2\beta_n R}-i\beta_n u_n w_{-n}e^{i2\beta_n R}\right) \nonumber\\
		&=0. \label{eqn:ImpedanceEstimateWUs}
	\end{align}
	Here, we have used the $ -\alpha $-quasiperiodic and $ \alpha $-quasiperiodic Rayleigh expansions for $ w $ and $ u^s $, respectively.
	The second integral can be further calculated as
	\begin{align}
		\left|\int_{\Gamma_R}\left(\partial_{\nu}w u^i-\partial_{\nu}u^i w\right) \,ds\right|  
		&= \left|\int_{0}^{2\pi} \left(\sum_{n \in \mathbb{Z}}i \hat{\beta}_n w_n e^{i(\hat{\alpha}_n x_1+\hat{\beta}_n R)}\gamma  e^{i(\alpha x_1-\beta R)} \right.\right. \nonumber\\
		&\quad \left.\left.+i\beta \gamma  e^{i(\alpha x_1-\beta R)}\sum_{n \in \mathbb{Z}} w_n e^{i(\hat{\alpha}_n x_1+\hat{\beta}_n R)}\right)\,dx_1\right| \nonumber\\
		&= \left|\int_{0}^{2\pi} \sum_{n \in \mathbb{Z}}i w_n \gamma  e^{i(\hat{\alpha}_n x_1+\hat{\beta}_n R)} e^{i(\alpha x_1-\beta R)}(\hat{\beta}_n+\beta)\,dx_1\right| \nonumber\\ 
		&=\left|\int_{0}^{2\pi} \sum_{n \in \mathbb{Z}}i w_n \gamma  (\hat{\beta}_n+\beta)e^{i n x_1} \cdot e^{i(\hat{\beta}_n-\beta)R}\,dx_1\right| \nonumber\\
		&= |2\pi w_0 \gamma (\hat{\beta}_0+\beta)| \nonumber\\
		&= |4\pi k \cos{\theta}w_0 \gamma | .\label{eqn:ImpedanceEstimateWUi}
	\end{align}
	Substituting (\ref{eqn:ImpedanceEstimateWUs}), (\ref{eqn:ImpedanceEstimateWUi}) and Lemma \ref{eqn:ImpedanceWvsW0} into (\ref{eqn:impedanceEstimateU}) and Lemma \ref{LemmaW} into the estimate of $\|u\|_{L^2(\Omega_R)}$, we get by applying Young's inequality that 
	\begin{align*}
		\|u\|_{L^2(\Omega_R)}^2 
		&\leq |4\pi k \cos{\theta}w_0 \gamma |+\|\partial_{\nu}w\|_{L^2(\Gamma)}\|u\|_{L^2(\Gamma)}\\
		&\leq C_9 \| w\|_{L^2(\Omega_R)}|\gamma |+\|\partial_{\nu}w\|_{L^2(\Gamma)}\|u\|_{L^2(\Gamma)}\\
		&\leq C_9 \tilde{C}\|u\|_{L^2(\Omega_R)}|\gamma |+\tilde{C}\|u\|_{L^2(\Omega_R)}\|u\|_{L^2(\Gamma)}\\
		&\leq C_9\left(\epsilon_3 |\tilde{C}|^2\|u\|_{L^2(\Omega_R)} ^2 +\frac{1}{4\epsilon_3}|\gamma |^2\right)+ \frac{1}{4}\|u\|_{L^2(\Omega_R)}^2+|\tilde{C}|^2\|u\|_{L^2(\Gamma)}^2,
	\end{align*}
	where $ \epsilon_3>0 $ is arbitrary and $  C_9=2\sqrt{2\pi}k \cos{\theta}$.
	Taking $\epsilon_3=\frac{1}{4 C_9|\tilde{C}|^2}$, we obtain
	\begin{align*}
		\|u\|_{L^2(\Omega_R)}^2 
		&\leq \frac{1}{4}\|u\|_{L^2(\Omega_R)}^2+C_9^2 |\tilde{C}|^2|\gamma |^2+ \frac{1}{4}\|u\|_{L^2(\Omega_R)}^2+|\tilde{C}|^2\|u\|_{L^2(\Gamma)}^2\\
		&= \frac{1}{2}\|u\|_{L^2(\Omega_R)}^2+|\tilde{C}|^2\|u\|_{L^2(\Gamma)}^2+C_9^2 |\tilde{C}|^2|\gamma |^2.
	\end{align*}
	Rearranging similar terms and using (\ref{eqn:estimatUgamma}), we have
	\begin{eqnarray*}
		\frac{1}{2}\|u\|_{L^2(\Omega_R)}^2 
		\leq |\tilde{C}|^2\|u\|_{L^2(\Gamma)}^2+C_9^2 |\tilde{C}|^2|\gamma |^2
		\leq \frac{|\tilde{C}|^2}{\lambda}|F(u)|+C_9^2 |\tilde{C}|^2|\gamma |^2.
	\end{eqnarray*}
	that is,
	\begin{eqnarray}
		\|u\|_{L^2(\Omega_R)}^2 \leq \frac{2|\tilde{C}|^2}{\lambda}|F(u)|+2C_9^2 |\tilde{C}|^2|\gamma |^2. \label{eqn:impedanceLlastU}
	\end{eqnarray}
	In view of (\ref{eqn:ImpedanceNabla}) and (\ref{eqn:impedanceLlastU}), we obtain
	\begin{align}
		\|u\|_{H^1_{\alpha}(\Omega_R)}^2&=k^2\|u\|_{L^2(\Omega_R)}^2+\|\nabla u\|_{L^2(\Omega_R)}^2 \nonumber\\
		&\leq 2k^2\|u\|_{L^2(\Omega_R)}^2+|F(u)| \nonumber\\
		&\leq 2k^2\left(\frac{2|\tilde{C}|^2}{\lambda}|F(u)|+2C_9^2 |\tilde{C}|^2|\gamma |^2\right)+|F(u)| \nonumber\\
		&= \left( 1+\frac{4k^2|\tilde{C}|^2}{\lambda}\right)|F(u)|+4k^2C_9^2 |\tilde{C}|^2|\gamma |^2. \label{eqn:ImpedanceUX1}
	\end{align}
	Similar to the Dirichlet case, we can also obtain (\ref{eqn:DirichletFuVsU0}), that is
	$$|F(u)|=|4\pi\beta \gamma \tilde{u}_0|.$$
	Recalling $\beta=k\cos{\theta}$, using Cauchy-Schwarz inequility and again using (\ref{eqn:U0UX}) (which in this case is $\|u\|_{H^1_{\alpha}(\Omega_R)}\geq  \sqrt{2\pi k}|\tilde{u}_0|  $), we deduce from (\ref{eqn:ImpedanceUX1}) that
	\begin{align}
		\|u\|_{H^1_{\alpha}(\Omega_R)}^2 
		&\leq \left( 1+\frac{4k^2|\tilde{C}|^2}{\lambda}\right)4\pi|k\cos{\theta}|| \gamma ||\tilde{u}_0|+4k^2C_9^2 |\tilde{C}|^2|\gamma |^2  \nonumber\\
		&= C_{10}|\gamma ||\tilde{u}_0|+C_{11}|\gamma |^2 \nonumber\\
		&\leq \frac{C_{10}\epsilon_4\|u\|_{H^1_{\alpha}(\Omega_R)}^2}{2\pi k} +\frac{C_{10}}{4\epsilon_4}|\gamma |^2+C_{11}|\gamma |^2, \label{eqn:ImpedanceUX2}
	\end{align}
	where 
	\begin{align*}
		C_{10}&=\left(1+\frac{4k^2|\tilde{C}|^2}{\lambda}\right)4\pi|k\cos{\theta}|,\\ 
		C_{11}&=4k^2C_9^2 |\tilde{C}|^2 =4k^2|\tilde{C}|^2 8\pi k^2\cos^2\theta=32\pi k^4\cos^2\theta |\tilde{C}|^2,
	\end{align*}
	and $ \epsilon_4>0 $ is arbitrary.
	Rearranging similar terms, we continue to calculate ($\ref{eqn:ImpedanceUX2}$) by
	\begin{eqnarray*}
		\left(1-\frac{C_{10}\epsilon_4}{2\pi k}\right)\|u\|_{H^1_{\alpha}(\Omega_R)}^2\leq\left(\frac{C_{10}}{4\epsilon_4}+C_{11}\right)|\gamma |^2.
	\end{eqnarray*}
	By taking $\epsilon_4=\frac{\pi k}{C_{10}}$, we obtain
	\begin{eqnarray*}
		\|u\|_{H^1_{\alpha}(\Omega_R)}^2\leq  \left(\frac{C_{10}}{2\epsilon_4}+2C_{11}\right)|\gamma |^2=\left(\frac{C_{10}^2}{2\pi k}+2C_{11}\right)|\gamma |^2.
	\end{eqnarray*}
	Furthermore, 
	\begin{align*}
		\frac{C_{10}^2}{2\pi k}+2C_{11}
		&= \left(1+\frac{4k^2|\tilde{C}|^2}{\lambda}\right)^2 8\pi k\cos^2\theta+64\pi k^4\cos^2\theta |\tilde{C}|^2\\
		&= 8\pi k\cos^2\theta\left[  \left( 1+\frac{4k^2|\tilde{C}|^2}{\lambda}\right)^2+8k^3|\tilde{C}|^2\right].
	\end{align*}
	Therefore,
	\ben
	\|u\|_{H^1_{\alpha}(\Omega_R)}^2\leq  8\pi k\cos^2\theta |\gamma |^2 C^{*2},
	\enn
	with
	\begin{eqnarray*}
		C^*=\sqrt{  \left( 1+\frac{4k^2\tilde{C}^2}{\lambda}\right)^2+8k^3|\tilde{C}|^2}.
	\end{eqnarray*}
	$\hfill\Box$
	
	\section{Transmission conditions} \label{sec:TC}
	The direct transmission scattering problem is to find the total $ u=u(x_1,x_2) $ such that 
	\begin{align}  
		\begin{array}{l}
			(TC):\qquad\left\{\begin{array}{ll}
				\Delta u+k^2u = 0    &\mbox{in} \quad \Omega_\Gamma^+, \\
				\Delta u+k_1^2u = 0    &\mbox{in} \quad \Omega_\Gamma^-, \\
				u^+=u^-,\;\partial_{\nu} u^+=\lambda\partial_{\nu} u^-          &\mbox{on} \quad \Gamma,\\
				u =u^i +u^s & \mbox{in} \quad \Omega_\Gamma^+,
			\end{array}\right. \\ \label{transmission}
		\end{array} 
	\end{align}
	here $ \lambda>0 $ is a constant and $ \nu $ is the normal direction at $ \Gamma $ pointing to $ \Omega_\Gamma^+ $. 
	We suppose that $ k,\;k_1 $ are constants and $ k\neq k_1 $. 
	For notational convenience, we set $ k_+=k,\;k_-=k_1,\; k(x)=k_{\pm}\;\rm{in}\; \Omega_\Gamma^{\pm}$. 
	The Rayleigh expansion radiation conditions together with proper outgoing radiation conditions as $ x_2 \to \pm\infty $ for penetrable gratings can be formulated as:
	\begin{align*}
		u^s&=\sum_{n \in \mathbb{Z}}u_n^+ e^{i\alpha_n x_1+i\beta_n^+ x_2},\quad \text{for} \quad  x_2>\Gamma_{\max}^+:=\max_{(x_1,x_2)\in\Gamma} x_2,\\
		u&=\sum_{n \in \mathbb{Z}}u_n^- e^{i\alpha_n x_1-i\beta_n^- x_2},\quad \text{for} \quad x_2<\Gamma_{\min}^-:=\min_{(x_1,x_2)\in\Gamma} x_2,
	\end{align*}
	where $ \alpha:=n+\alpha $ and
	\begin{align*}
		\begin{array}{l}
			\beta_n^+=\left\{\begin{array}{ll}
				\sqrt{k_+^2-|\alpha_n|^2}, & |\alpha_n|\leq k_+, \\
				i\sqrt{|\alpha_n|^2-k_+^2}, &  |\alpha_n|> k_+,
			\end{array}\right. \\ 
		\end{array}
		\begin{array}{l}
			\beta_n^-=\left\{\begin{array}{ll}
				\sqrt{k_-^2-|\alpha_n|^2}, & |\alpha_n|\leq k_-, \\
				i\sqrt{|\alpha_n|^2-k_-^2}, &  |\alpha_n|> k_-.
			\end{array}\right. \\ 
		\end{array} 
	\end{align*}
	Throughout this section, we suppose that $ f_-<\min_{x\in\Gamma} x_2<\max_{x\in\Gamma} x_2< f_+$, where $ f_-<f_+ $.
	Then we introduce the notations
	\begin{align*}
		S_R=S_{R,f}:=\{x\in\mathbb{R}^2:|f(x_1)|<|x_2|<R,0<x_1<2\pi\} \quad\text{for some }  R>\max\{|f_+|,|f_-|\},	
	\end{align*}
	$ S_R^{\pm}:=\{x\in \Omega_\Gamma^{\pm}: -R<x_2<R\} $ and $ \Gamma_{R}^{\pm}:=\{(x_1,\pm R):0<x_1<2\pi\} $. 
	Define the DtN mappings $ T^{\pm} :H^{1/2}_{\alpha}(\Gamma_{R}^{\pm}) \to H^{-1/2}_{\alpha}(\Gamma_{R}^{\pm}) $ as 
	\begin{align*}
		T^{\pm}g(x_1):=\pm\sum_{n \in \mathbb{Z}}i\beta^{\pm}_n g_n e^{i\alpha_n x_1}, \qquad g(x_1)=\sum_{n \in \mathbb{Z}} g_n e^{i\alpha_n x_1}\in H_{\alpha}^{1/2}(\Gamma^{\pm}_R).
	\end{align*}
	We also define a piecewise constant function $ a(x) $ such that $ a(x)=1 $	in $ S_R^+ $, $ a(x)=\lambda $ in $ S_R^- $, therefore 
	\ben
	\|\partial_2 u\|^2_{L^2(S_R)}=\int_{S_R}a(x)|\partial_2 u|^2\,dx,\quad \| u\|^2_{L^2(S_R)}=\int_{S_R}a(x)|u|^2\,dx.
	\enn
	Then the Sobolev space $ H^1_\alpha(S_R) $ is equipped with the norm 
	\be
	\| u\|^2_{H^1_\alpha(S_R)}:=\int_{S_R} a(x)|\nabla u|^2+a(x) k^2(x)|u|^2\,dx. \label{trNorm}
	\en
	From \cite{xiangJianli2021} we know that problem $ (\ref{transmission}) $ is uniquely solvable in $ H^1_\alpha(S_R) $ and then we give the following theorem.
	\begin{theorem}\label{TransTH}
		Let $u^i$ be a plane wave given by (\ref{eqn:Ui}) and suppose that $u \in H^1_\alpha(S_R)$ be the unique solution to the Dirichlet problem (\ref{transmission}). 
		Choose $ f_+,f_-\in\mathbb{R} $ such that $ \min_{x_1\in (0,2\pi)}f(x_1)-f_->1 $ and $ \max_{x_1\in (0,2\pi)}f(x_1)-f_+<-1 $.\\
		$ \text{(i).} $ If $ \lambda \geq 1,\; k^2_+>\lambda  k^2_- $, we have
		$$\|u\|_{H^1_\alpha(S_R)} \leq 2\sqrt{2\pi k} \cos\theta C_{12} |\gamma |,$$
		where 
		\begin{align*}
			C_{12}&=2\max\{k_+,\lambda  k_-\}\frac{[2k_+(R-f_-)+1]}{C_T}+1,\\
			C_T &= \min \left \{\frac{2}{( R-f_-)^2},\frac{(k_+^2-\lambda k_-^2)}{2( R-f_-)\sqrt{1+L^2}}\right \}.
		\end{align*}
		$ \text{(ii)}. $ If $ \lambda \leq 1,\; k^2_+<\lambda  k^2_- $, we have 
		$$\|u\|_{H^1_\alpha(S_R)} \leq 2\sqrt{2\pi k}\cos\theta C_{13} |\gamma |,$$
		where 
		\begin{align*}
			C_{13}&=2\max\{k_+,\lambda k_-\}\frac{[2k_+(R-f_+)+1]}{C_S}+1, \\
			C_S&=\min \left\{\frac{2}{( R-f_-)^2},\frac{(\lambda k_-^2-k_+^2)}{2( R-f_-)\sqrt{1+L^2}}\right\} ,	
		\end{align*}
		Here $L$ is the Lipschitz constant of $f$.
	\end{theorem}
	In order to prove the theorem, we need some preliminaries. 
	We first deduce the variational formulation.
	Multiplying on both sides of the Helmholtz equation by the complex conjugate of $ v $ and integrating over $ S_R^+ $ and $ S_R^- $ respectively, and using the boundary condition, we obtain
	\begin{align}
		\int_{S_R^+} \nabla u \cdot \nabla \Bar{v}-k_+^2u\Bar{v}\, dx-\int_{\Gamma_R^+}(T^+u)\Bar{v}\, ds+\int_{\Gamma}\partial_\nu u^+\bar{v}\, ds&=-\int_{\Gamma_R^+}(T^+u^i-\partial_2 u^i)\Bar{v}\, ds,\label{Imvf1}\\
		\int_{S_R^-} \nabla u \cdot \nabla \Bar{v}-k_-^2u\Bar{v}\, dx+\int_{\Gamma_R^-}(T^-u)\Bar{v}\, ds-\int_{\Gamma}\partial_\nu u^-\bar{v}\, ds&=0.\label{Imvf2}
	\end{align}
	Multiplying (\ref{Imvf2}) by $ \lambda $ and adding it to (\ref{Imvf1}), we have
	\begin{align*}
		&\;\int_{S_R^+} \nabla u \cdot \nabla \Bar{v}-k_+^2u\Bar{v}\, dx+\lambda \int_{S_R^-} \nabla u \cdot \nabla \Bar{v}-k_-^2u\Bar{v}\, dx-\int_{\Gamma_R^+}(T^+u)\Bar{v}\, ds \\
		&\quad\quad+\lambda \int_{\Gamma_R^-}(T^-u)\Bar{v}\, ds
		=-\int_{\Gamma_R^+}(T^+u^i-\partial_2 u^i)\Bar{v}\, ds,
	\end{align*}
	where we have used the transmission condition $ \partial_\nu u^+ -\lambda\partial_\nu u^-=0 $ on $ \Gamma $. 
	Straightforward calculations show that
	\ben
	\int_{\Gamma_R^+}(T^+u^i-\partial_2 u^i)\Bar{v}\, ds=2i\beta e^{-i\beta R}\gamma \int_0^{2\pi}e^{i\alpha x_1}\Bar{v}(x_1,R)\, dx_1.
	\enn
	Besides, by the definition of $ k (x)$ and $ a(x) $, we have
	\begin{align*}
		\int_{S_R^+} \nabla u \cdot \nabla \Bar{v}-k_+^2u\Bar{v}\, dx+\lambda \int_{S_R^-} \nabla u \cdot \nabla \Bar{v}-k_-^2u\Bar{v}\, dx=\int_{S_R}a(x) \nabla u \cdot \nabla \Bar{v}-a(x)k^2(x)u\Bar{v}\, dx.
	\end{align*}
	Therefore the variational formulation for the scattering problem can be written as: find $ u\in H^1_\alpha (S_R) $ such that 
	\begin{equation}
		a(u,v)=F(v)  \quad \text{ for all }  v \in H^1_{\alpha}(S_R),   \label{eqn:varitional formulaTRans}
	\end{equation}
	where
	\begin{align*}
		a(u,v)&:=\int_{S_R}a(x) \nabla u \cdot \nabla \Bar{v}-a(x)k^2(x)u\Bar{v}\, dx-\int_{\Gamma_R^+}(T^+u)\Bar{v}\, ds+\lambda\int_{\Gamma_R^-}(T^-u)\Bar{v}\, ds, \\
		F(v)&:=-2i\beta e^{-i\beta R}\gamma \int_0^{2\pi}e^{i\alpha x_1}\Bar{v}(x_1,R)\, dx_1. 
	\end{align*}	
	Taking the real and imaginary parts on both sides of $ (\ref{eqn:varitional formulaTRans}) $	with $ v=u $, we get
	\be
	\int_{S_R}\left[a(x)|\nabla u|^2-a(x)k^2(x)|u|^2\right]\,dx-\mbox{Re}\int_{\Gamma_R^+}(T^+u)\Bar{u}\, ds+\lambda\mbox{Re}\int_{\Gamma_R^-}(T^-u)\Bar{u}\, ds=\mbox{Re}\,F(u),\label{TranRe}\\
	-\mbox{Im}\int_{\Gamma_R^+}(T^+u)\Bar{u}\, ds+\lambda\mbox{Im}\int_{\Gamma_R^-}(T^-u)\Bar{u}\, ds=\mbox{Im}\,F(u).\label{TranIm}
	\en
	By definition of the DtN maps $ T^{\pm} $, we can easily get for $ u\big|_{\Gamma_R^{\pm}}=\sum_{n\in \mathbb{Z}}\hat{u}_n^\pm e^{i\alpha_n x_1} $ that
	\begin{align}
		\mbox{Re}\,\int_{\Gamma_R^{\pm}}\pm(T^{\pm}u)\Bar{u}&=\mbox{Re}\,2\pi \sum_{n\in \mathbb{Z}}i\beta_n^{\pm}|\hat{u}_n^\pm|^2=-2\pi \sum_{|\alpha_n|>k_\pm}|\beta_n^{\pm}||\hat{u}_n^\pm|^2 \leq 0,\label{trRe}\\
		\mbox{Im}\,\int_{\Gamma_R^{\pm}}\pm(T^{\pm}u)\Bar{u}&=\mbox{Im}\,2\pi \sum_{n\in \mathbb{Z}}i\beta_n^{\pm}|\hat{u}_n^\pm|^2=2\pi \sum_{|\alpha_n|\leq k_\pm}|\beta_n^{\pm}||\hat{u}_n^\pm|^2 \geq 0.\label{trIm}
	\end{align}
	We also need the following lemma.
	\begin{lemma}\label{RellichTrans}
		If $u\in H_\alpha^1(S_R)$ is a solution to (\ref{transmission}), then we have
		\begin{align*}
			&\left( \int_{\Gamma_R^+}- \int_{\Gamma}\right) (x_2-c)\left[-\nu_2|\nabla u^+|^2+\nu_2k^2_+|u|^2+2\mbox{Re}\,(\partial_2\bar{u}^+\partial_\nu u^+)\right]\,ds\\
			&+\lambda \left( \int_{\Gamma_R^-}+ \int_{\Gamma}\right) (x_2-c)\left[-\nu_2|\nabla u^-|^2+\nu_2k^2_-|u|^2+2\mbox{Re}\,(\partial_2\bar{u}^-\partial_\nu u^-)\right]\,ds\\
			&+\int_{S_R^+}|\nabla u|^2-k^2_+|u|^2-2|\partial_2 u|^2\,dx+\lambda \int_{S_R^-}|\nabla u|^2-k^2_-|u|^2-2|\partial_2 u|^2\,dx=0,
		\end{align*}
		where the normal direction $ \nu=(\nu_1,\nu_2) \in \mathbb{S}:=\{x\in \mathbb{R}^2: |x|=1\} $ at $ \Gamma $ is supposed to point to $ S_R^+ $ and $ c $ is a constant.
	\end{lemma}
	\begin{proof}
		Set $ \Gamma_{j}=\{x_2=f_j(x_1),\,x_1\in (0,2\pi)\} $.
		Following Remark 3.2 discussed in \cite{Elschner2002}, there exists a sequence of infinitely smooth profile functions $ f_j,j\in \mathbb{N} $ such that $ \|f_j-f\|_{C_p^{0,1}}\to\infty $ as $ j \to\infty$.
		Denote by $ u_j \in H_\alpha^1(S_R) $ the unique solution of transmission problem (\ref{transmission}) with $ S_{R,j}^\pm:=S_{R,f_j}^\pm $.
		Then we establish the equation in this lemma when $ \Gamma=\Gamma_{f_j} $ is an infinitely smooth profile. 
		For $ u_j\in  H_\alpha^2(S_{R,j}^+) \cap H_\alpha^2(S_{R,j}^-) $ , using the Rellich identity (see Lemma \ref{eqn:rellichNohelmholtz} in the Appendix and note here that $ u_j $ satisfies the Helmholtz equation) on both sides of $ S_{R,j}^+ $ and $ S_{R,j}^- $, noting that $ \nu_2 \geq 0 $, we obtain
		\begin{align*}
			0&=\left( \int_{\Gamma_R^{\pm}}\mp \int_{\Gamma_j}\right) (x_2-c)\left[-\nu_2|\nabla u_j^\pm|^2+\nu_2k^2_\pm|u_j|^2+2\mbox{Re}\,(\partial_2\bar{u}_j^\pm\partial_\nu u_j^\pm)\right]\,ds\\
			&\quad+\int_{S_R^\pm}|\nabla u_j|^2-k^2_\pm|u_j|^2-2|\partial_2 u_j|^2\,dx :=I_c^{\pm}.
		\end{align*}
		Therefore, 
		\begin{align}
			0 &= I_c^++\lambda I_c^- \nonumber\\
			&=\left( \int_{\Gamma_R^+}- \int_{\Gamma_j}\right) (x_2-c)\left[-\nu_2|\nabla u_j^+|^2+\nu_2k^2_+|u_j|^2+2\mbox{Re}\,(\partial_2\bar{u}_j^+\partial_\nu u_j^+)\right]\,ds \nonumber\\
			&\quad+\lambda \left( \int_{\Gamma_R^-}+ \int_{\Gamma_j}\right) (x_2-c)\left[-\nu_2|\nabla u_j^-|^2+\nu_2k^2_-|u_j|^2+2\mbox{Re}\,(\partial_2\bar{u}_j^-\partial_\nu u_j^-)\right]\,ds \nonumber\\
			&\quad+\int_{S_R^+}|\nabla u_j|^2-k^2_+|u_j|^2-2|\partial_2 u_j|^2\,dx+\lambda \int_{S_R^-}|\nabla u_j|^2-k^2_-|u_j|^2-2|\partial_2 u_j|^2\,dx. \label{ICRellichJ}
		\end{align}
		Hence passing the limit $ j \to\infty$ in the inequality (\ref{ICRellichJ}) and using the convergence $  u_j \to u$ in $ H_\alpha^1(S_R) $ we get a similar equation for $ u  $. 
		Using the same trick just as the end of the proof in Lemma $ \ref{LemmaW} $, we obtain the identity in the lemma, which finishes the proof.
	\end{proof}
	\textbf{Proof of Theorem \ref{TransTH}.}
	\text{(i).}	Suppose $ \lambda \geq 1,\; k^2_+>\lambda  k^2_- $. 
	Straightforward calculations show that for $ u\big|_{\Gamma_R^{\pm}}=\sum_{n\in \mathbb{Z}}\hat{u}_n^\pm e^{i\alpha_n x_1} $, we obtain
	\begin{align}
		&\quad \int_{\Gamma_R^{\pm}}  (x_2-f_-)\left[-\nu_2|\nabla u^\pm|^2+\nu_2k^2_\pm|u|^2+2\mbox{Re}\,(\partial_2\bar{u}^\pm\partial_\nu u^\pm)\right]\,ds \no\\
		&=\int_{\Gamma_R^{\pm}}(R-f_-)\left[|\partial_2u^\pm|^2-|\partial_1u^\pm|^2+k^2_\pm|u|^2\right]\,ds \no\\
		&=(\pm R-f_-)4\pi \sum_{|\alpha_n|\leq k_\pm}|\beta_n^{\pm}|^2|u_n^\pm|^2. \label{REsti}
	\end{align}
	Simple calculations show $ \mbox{Re}\,(\partial_2\bar{u}^\pm\partial_\nu u^\pm)=\nu_2|\partial_\nu u^\pm|^2 $.
	By the transmission condition, we have
	\begin{align}
		&\quad-|\nabla u^+|^2+\lambda |\nabla u^-|^2+2|\partial_\nu  u^+|^2-2\lambda |\partial_\nu  u^-|^2 \no \\
		&=-\left(|\partial_\nu  u^+|^2+|\partial_\tau  u^+|^2\right)+\lambda \left(|\partial_\nu  u^-|^2+|\partial_\tau  u^-|^2\right)+2|\partial_\nu  u^+|^2-2\lambda |\partial_\nu  u^-|^2 \no \\
		&=-\left(\lambda^2 |\partial_\nu  u^-|^2+|\partial_\tau  u^-|^2\right)+\lambda\left(|\partial_\nu  u^-|^2+|\partial_\tau  u^-|^2\right)+2\lambda^2 |\partial_\nu  u^-|^2-2\lambda |\partial_\nu  u^-|^2  \no \\
		&=(\lambda^2-\lambda)|\partial_\nu  u^-|^2+(\lambda-1)|\partial_\tau  u^-|^2 \label{relation}
	\end{align}
	Applying Lemma $ \ref{RellichTrans} $ with $ c=f_- $ and (\ref{REsti}), (\ref{relation}), we have
	\begin{align}
		0 &=\left( \int_{\Gamma_R^+}- \int_{\Gamma}\right) (x_2-f_-)\left[-\nu_2|\nabla u^+|^2+\nu_2k^2_+|u|^2+2\mbox{Re}\,(\partial_2\bar{u}^+\partial_\nu u^+)\right]\,ds \no\\
		&\quad+\lambda \left( \int_{\Gamma_R^-}+ \int_{\Gamma}\right) (x_2-f_-)\left[-\nu_2|\nabla u^-|^2+\nu_2k^2_-|u|^2+2\mbox{Re}\,(\partial_2\bar{u}^-\partial_\nu u^-)\right]\,ds \no\\
		&\quad+\int_{S_R^+}|\nabla u|^2-k^2_+|u|^2-2|\partial_2 u|^2\,dx+\lambda \int_{S_R^-}|\nabla u|^2-k^2_-|u|^2-2|\partial_2 u|^2\,dx \no\\
		&=( R-f_-)4\pi \sum_{|\alpha_n|\leq k_+}|\beta_n^+|^2|u_n^+|^2+\lambda (- R-f_-)4\pi \sum_{|\alpha_n|\leq k_-}|\beta_n^-|^2|u_n^-|^2 \no\\
		&\quad- \int_{\Gamma}\left[\lambda(\lambda-1)|\partial_\nu u^-|^2+(\lambda-1)|\partial_\tau u^-|^2+(k_+^2-\lambda k_-^2)|u|^2\right]\nu_2 (x_2-f_-)\,ds \no\\
		&\quad-2\int_{S_R}a(x)|\partial_2 u|^2\,dx+\int_{S_R}a(x)|\nabla u|^2-a(x)k^2(x)|u|^2\,dx. \label{fstEsti}
	\end{align}
	Using $ (\ref{TranRe}) $ and (\ref{fstEsti}), we have
	\begin{align}
		&\quad\int_{\Gamma}\left[\lambda(\lambda-1)|\partial_\nu u^-|^2+(\lambda-1)|\partial_\tau u^-|^2+(k_+^2-\lambda k_-^2)|u|^2\right]\nu_2 (x_2-f_-)\,ds+2\int_{S_R}a(x)|\partial_2 u|^2\,dx \nonumber\\
		&=( R-f_-)4\pi \sum_{|\alpha_n|\leq k_+}|\beta_n^+|^2|\hat{u}_n^+|^2+\lambda (- R-f_-)4\pi \sum_{|\alpha_n|\leq k_-}|\beta_n^-|^2|\hat{u}_n^-|^2 \nonumber\\
		&\quad+\int_{S_R}a(x)|\nabla u|^2-a(x)k^2(x)|u|^2\,dx \nonumber\\
		&=4\pi( R-f_-)\sum_{|\alpha_n|\leq k_+}|\beta_n^+|^2|\hat{u}_n^+|^2+4\pi\lambda (- R-f_-)\sum_{|\alpha_n|\leq k_-}|\beta_n^-|^2|\hat{u}_n^-|^2 \nonumber\\
		&\quad+ \mbox{Re}\,\int_{\Gamma_R^+}(T^+u)\Bar{u}\,ds-\lambda \mbox{Re}\,\int_{\Gamma_R^-}(T^-u)\Bar{u}\,ds +\mbox{Re}\, F(u).\label{TranEsti1}
	\end{align}
	From $ (\ref{TranIm}) $, $ (\ref{trRe}) $ and $ (\ref{trIm}) $, we obtain
	\begin{align}
		\mbox{Im}\, F(u)&=-\mbox{Im}\,\int_{\Gamma_R^+}(T^+u)\Bar{u}\,ds+\lambda \mbox{Im}\,\int_{\Gamma_R^-}(T^-u)\Bar{u}\,ds\no\\
		&=-2\pi\sum_{|\alpha_n|\leq k_+}|\beta_n^{+}||\hat{u}_n^+|^2-2\pi\lambda\sum_{|\alpha_n|\leq k_-}|\beta_n^{-}||\hat{u}_n^-|^2\leq 0, \label{trFu}
	\end{align}
	and
	\begin{align}
		\mbox{Re}\,\int_{\Gamma_R^+}(T^+u)\Bar{u}\,ds-\lambda \mbox{Re}\,\int_{\Gamma_R^-}(T^-u)\Bar{u}\,ds \leq 0. \label{rtr}
	\end{align}
	We also need the following estimate. 
	Using (\ref{trFu}), $|\beta_n^+|\leq k_+ $ and $ R\pm f_->0 $, we have
	\begin{align*}
		&\quad 4\pi( R-f_-)\sum_{|\alpha_n|\leq k_+}|\beta_n^+|^2|\hat{u}_n^+|^2+4\pi\lambda (- R-f_-)\sum_{|\alpha_n|\leq k_-}|\beta_n^-|^2|\hat{u}_n^-|^2+2k_+(R-f_-)\mbox{Im}\, F(u)\no\\
		&= 4\pi \sum_{|\alpha_n|\leq k_+}( R-f_-)(|\beta_n^+|-k_+)|\beta_n^+| |\hat{u}_n^+|^2 \no\\
		&\quad +4\pi\lambda \sum_{|\alpha_n|\leq k_-}\left[(- R-f_-)|\beta_n^-|-k_+(R-f_-)\right]|\beta_n^-||\hat{u}_n^-|^2 \leq 0.
	\end{align*}
	Recalling that the conditions $ \lambda \geq 1 $ and $ k_+^2> \lambda k_-^2 $, $ \nu_2\geq C_L $ and $ \min_{x_1\in (0,2\pi)}f(x_1)-f_->1 $, we can get a low bound of the left-hand side of $ (\ref{TranEsti1}) $,
	\begin{align}
		&\int_{\Gamma}\left[\lambda(\lambda-1)|\partial_\nu u^-|^2+(\lambda-1)|\partial_\tau u^-|^2+(k_+^2-\lambda k_-^2)|u|^2\right]\nu_2 (x_2-f_-)\,ds \no\\
		&\quad +2\int_{S_R}a(x)|\partial_2 u|^2\,dx
		\geq (k_+^2-\lambda k_-^2)C_L\|u\|^2_{L^2(\Gamma)}+2\|\partial_2 u\|^2_{L^2(S_R)}, \label{lEsti1}
	\end{align}
	where $ C_L=\frac{1}{\sqrt{1+L^2}} $.
	Therefore, combining $ (\ref{TranEsti1}) $ and $ (\ref{rtr}) $-$ (\ref{lEsti1}) $ yields
	\begin{align}
		(k_+^2-\lambda k_-^2)C_L\|u\|^2_{L^2(\Gamma)}+2\|\partial_2 u\|^2_{L^2(S_R)} &\leq -2k_+(R-f_-)\mbox{Im}\, F(u)+\mbox{Re}\, F(u)\no\\
		&\leq[2k_+(R-f_-)+1]|F(u)|. \label{trEs}
	\end{align}
	Using Lemma \ref{eqn:estimator2}, the $ L^2 $-norm of $ u $ can be bounded by $(\ref{trEs})  $ as follows:
	\begin{align*}
		C_T\|u\|^2_{L^2(S_R)} &\leq C_T\left[( R-f_-)^2\|\partial_2 u\|^2_{L^2(S_R)}+2( R-f_-)\|u\|^2_{L^2(\Gamma)}\right] \\
		&\leq 2\|\partial_2 u\|^2_{L^2(S_R)}+(k_+^2-\lambda k_-^2)C_L\|u\|^2_{L^2(\Gamma)},
	\end{align*}
	where $ C_T=\min \left \{\frac{2}{( R-f_-)^2},\frac{(k_+^2-\lambda k_-^2)C_L}{2( R-f_-)}\right \} $.
	Therefore, using $ (\ref{trEs}) $,
	\ben
	C_T\|u\|^2_{L^2(S_R)}\leq [2k_+(R-f_-)+1]|F(u)|,
	\enn
	from which we obtain
	\ben
	\|u\|^2_{L^2(S_R)}\leq \frac{[2k_+(R-f_-)+1]}{C_T}|F(u)|.
	\enn
	Then by (\ref{trNorm}), (\ref{TranRe}) and (\ref{trFu}), we have
	\begin{align}
		\| u\|^2_{H^1_\alpha(S_R)}&=\int_{S_R} a(x)|\nabla u|^2+a(x) k^2(x)|u|^2\,dx \no\\
		&=\int_{S_R}2a(x)k^2(x)|u|^2\,dx+\mbox{Re}\,\int_{\Gamma_R^+}(T^+u)\Bar{u}\, ds-\lambda\mbox{Re}\,\int_{\Gamma_R^-}(T^-u)\Bar{u}\, ds+\mbox{Re}\,F(u) \no\\
		&\leq  C_{12}|F(u)| \label{fstUx}
	\end{align}
	where $ C_{12}=2\max\{k_+,\lambda k_-\}\frac{[2k_+(R-f_-)+1]}{C_T}+1 $. 
	As done in the Dirichlet case, we can estimate $ |F(u)| $ by (see (\ref{eqn:DirichletFuVsU0}) and (\ref{eqn:U0UX}))
	\ben
	|F(u)|=4\pi \beta |\gamma  \tilde{u}_0|\leq 4\pi \beta \left( \epsilon_5 |\tilde{u}_0|^2+\frac{1}{4\epsilon_5}|\gamma |^2\right) \leq 2 \epsilon_5\cos\theta\| u\|^2_{H^1_\alpha(S_R)}+\frac{\pi \beta}{\epsilon_5}|\gamma |^2.
	\enn
	Therefore,
	\ben
	\| u\|^2_{H^1_\alpha(S_R)}\leq C_{12}\left(2 \epsilon_5\cos\theta\| u\|^2_{H^1_\alpha(S_R)}+\frac{\pi \beta}{\epsilon_5}|\gamma |^2 \right).
	\enn
	Taking $ \epsilon_5=\frac{1}{4 \cos\theta C_{12}} $, we have
	\ben
	\| u\|^2_{H^1_\alpha(S_R)}\leq 8\pi k \cos^2 \theta  C_{12}^2  |\gamma |^2.
	\enn
	\text{(ii).}Suppose that $ \lambda \leq 1$ and $ k^2_+<\lambda  k^2_- $. 
	In Lemma $ \ref{RellichTrans} $, we take $ c=f_+ $. 
	Hence, similar to (\ref{REsti}) and (\ref{fstEsti}), we get
	\ben
	\int_{\Gamma_R^{\pm}}  (x_2-f_+)\left[-\nu_2|\nabla u^\pm|^2+\nu_2k^2_\pm|u|^2+2\mbox{Re}\,(\partial_2\bar{u}^\pm\partial_\nu u^\pm)\right]\,ds=(\pm R-f_+)4\pi \sum_{|\alpha_n|\leq k_\pm}|\beta_n^{\pm}|^2|\hat{u}_n^\pm|^2,
	\enn
	and
	\begin{align*}
		0 &=4\pi( R-f_+)\sum_{|\alpha_n|\leq k_+}|\beta_n^+|^2|\hat{u}_n^+|^2+4\pi\lambda (- R-f_+)\sum_{|\alpha_n|\leq k_-}|\beta_n^-|^2|\hat{u}_n^-|^2\\
		&\quad-2\int_{S_R}a(x)|\partial_2 u|^2\,dx+\int_{S_R}a(x)|\nabla u|^2-a(x)k^2(x)|u|^2\,dx\\
		&\quad - \int_{\Gamma}\left[\lambda(\lambda-1)|\partial_\nu u^-|^2+(\lambda-1)|\partial_\tau u^-|^2+(k_+^2-\lambda k_-^2)|u|^2\right]\nu_2 (x_2-f_+)\,ds.\\
	\end{align*}
	Therefore, similar to (\ref{TranEsti1}), we get
	\begin{align}
		&\quad\int_{\Gamma}\left[\lambda(\lambda-1)|\partial_\nu u^-|^2+(\lambda-1)|\partial_\tau u^-|^2+(k_+^2-\lambda k_-^2)|u|^2\right]\nu_2 (x_2-f_+)\,ds+2\int_{S_R}a(x)|\partial_2 u|^2\,dx \nonumber\\
		&=4\pi( R-f_+)\sum_{|\alpha_n|\leq k_+}|\beta_n^+|^2|\hat{u}_n^+|^2+4\pi\lambda (- R-f_+)\sum_{|\alpha_n|\leq k_-}|\beta_n^-|^2|\hat{u}_n^-|^2 \no\\
		&\quad + \mbox{Re}\,\int_{\Gamma_R^+}(T^+u)\Bar{u}\,ds-\lambda \mbox{Re}\,\int_{\Gamma_R^-}(T^-u)\Bar{u}\,ds +\mbox{Re}\, F(u).\label{TranEsti2}
	\end{align}	
	Besides, by (\ref{trFu}),
	\begin{align}
		&\quad 4\pi( R-f_+)\sum_{|\alpha_n|\leq k_+}|\beta_n^+|^2|\hat{u}_n^+|^2+4\pi\lambda (- R-f_+)\sum_{|\alpha_n|\leq k_-}|\beta_n^-|^2|\hat{u}_n^-|^2
		+2k_+(R-f_+)\mbox{Im}\, F(u) \no\\
		&= 4\pi \sum_{|\alpha_n|\leq k_+}( R-f_+)(|\beta_n^+|-k_+)|\beta_n^+| |\hat{u}_n^+|^2 \no\\
		&\quad +4\pi\lambda \sum_{|\alpha_n|\leq k_-}\left[(- R-f_+)|\beta_n^-|-k_+(R-f_+)\right]|\beta_n^-||\hat{u}_n^-|^2 \leq 0.\label{scdEsti}
	\end{align}
	Therefore, one deduces from (\ref{TranEsti2}), (\ref{scdEsti}) and $ \max_{x_1\in (0,2\pi)}f(x_1)-f_+<-1 $ that
	\begin{align*}
		(\lambda k_-^2-k_+^2)C_L\|u\|^2_{L^2(\Gamma)}+2\|\partial_2 u\|^2_{L^2(S_R)} &\leq -2k_+( R-f_+)\mbox{Im}\, F(u)+\mbox{Re}\, F(u)\\
		&\leq [2k_+( R-f_+)+1]|F(u)|.
	\end{align*}
	We continue to estimate the left-hand side of the above equation, 
	\begin{align*}
		C_S\|u\|^2_{L^2(S_R)} &\leq C_S\left[( R-f_-)^2\|\partial_2 u\|^2_{L^2(S_R)}+2( R-f_-)\|u\|^2_{L^2(\Gamma)}\right] \\
		&\leq 2\|\partial_2 u\|^2_{L^2(S_R)}+(\lambda k_-^2-k_+^2)C_L\|u\|^2_{L^2(\Gamma)},
	\end{align*}
	where $ C_S=\min \left \{\frac{2}{( R-f_-)^2},\frac{(\lambda k_-^2-k_+^2)C_L}{2( R-f_-)}\right \} $.
	Therefore,
	\ben
	C_S\|u\|^2_{L^2(S_R)}\leq  [2k_+( R-f_+)+1]|F(u)|,
	\enn
	from which we obtain
	\ben
	\|u\|^2_{L^2(S_R)}\leq \frac{ [2k_+( R-f_+)+1]}{C_S}|F(u)|.
	\enn
	Then by $ (\ref{TranRe}) $, we have (see, e.g. (\ref{fstUx}))
	\ben
	\| u\|^2_{H^1_\alpha(S_R)}\leq 2\max \{k_+,\lambda k_-\}\|u\|^2_{L^2(S_R)}+|F(u)| \leq
	C_{13}|F(u)|
	\enn
	where $ C_{13}=2\max\{k_+,\lambda k_-\}\frac{ [2k_+( R-f_+)+1]}{C_S}+1 $. 
	As in the previous method, we then estimate $ |F(u)| $.
	\ben
	|F(u)|=4\pi \beta |\gamma  \tilde{u}_0|\leq 4\pi \beta \left( \epsilon_6 |\tilde{u}_0|^2+\frac{1}{4\epsilon_6}|\gamma |^2\right) \leq 2 \epsilon_6\cos\theta \| u\|^2_{H^1_\alpha(S_R)}+\frac{\pi \beta}{\epsilon_6}|\gamma |^2.
	\enn
	Therefore,
	\ben
	\| u\|^2_{H^1_\alpha(S_R)}\leq C_{13}\left(2 \epsilon_6\cos\theta \| u\|^2_{H^1_\alpha(S_R)}+\frac{\pi \beta}{\epsilon_6}|\gamma |^2\right),
	\enn
	by taking $ \epsilon_6=\frac{1}{4 \cos\theta C_{13}} $, we have
	\ben
	\| u\|^2_{H^1_\alpha(S_R)}\leq 8\pi k \cos^2 \theta  C_{13}^2  |\gamma |^2.
	\enn
	$\hfill\Box$	
	
	\section{Appendix} \label{sec:App}
	\begin{lemma}\label{eqn:rellichNohelmholtz}
		
		If $v \in H_{\alpha}^2(\Omega_R)$, then we have the Rellich identity 
		\ben
		2{\rm{Re}}\,\int_{\Omega_R}(x_2-c)\partial_2 \Bar{v} (\Delta v+k^2v)\,dx-\int_{\Omega_R}|\nabla v|^2-k^2|v|^2-2|\partial_2 v|^2\,dx\nonumber\\	
		=\left(\int_{\Gamma_R}-\int_{\Gamma}\right)(x_2-c)[-\nu_2|\nabla v|^2+\nu_2k^2|v|^2+2{\rm{Re}}\,(\partial_2\bar{v}\partial_\nu v)]\,ds,\qquad
		\enn
		where $ c $ is a constant and $ \nu=(\nu_1,\nu_2) $ is the normal direction at $ \Gamma \cup \Gamma_{R} $ pointing upward.
	\end{lemma}
	
	\begin{proof}
		Using the Green's formula and the trick $2 \mbox{Re}\,(\partial_2 \Bar{v} v)=\partial_2 |v|^2$ (note that $ \nu=(0,1)\;\rm{on}\;\Gamma_{R} $), we have
		\begin{align}
			&\quad 2\mbox{Re} \,\int_{\Omega_R} (x_2-c)\partial_2 \Bar{v} (\Delta v+k^2 v)\,dx \nonumber \\ 
			&= 2\mbox{Re}\,\int_{\Omega_R} -\nabla((x_2-c)\partial_2\Bar{v})\cdot \nabla v+k^2(x_2-c)\partial_2 \Bar{v} v \,dx +2\mbox{Re}\,\int_{\partial\Omega_R}(x_2-c)\partial_2\Bar{v} \partial_{\nu}v\,ds \nonumber\\
			&=-2\mbox{Re}\,\int_{\Omega_R} (0,\partial_2\Bar{v})\cdot (\partial_1 v,\partial_2 v)^\top+(x_2-c)\partial_2\nabla\Bar{v}\cdot\nabla v\,dx+\int_{\Omega_R} k^2 (x_2-c)\partial_2 |v|^2 \,dx\nonumber\\
			&\quad+2\mbox{Re}\,\int_{\Gamma_R}(x_2-c)\nu_2|\partial_\nu v|^2\,ds -2\mbox{Re}\,\int_{\Gamma}(x_2-c)\nu_2|\partial_{\nu}v|^2\,ds\nonumber\\
			&=-2\int_{\Omega_R }|\partial_2 v|^2\,dx-\int_{\Omega_R}(x_2-c)\partial_2|\nabla v|^2\,dx+\int_{\Omega_R}k^2(x_2-c)\partial_2|v|^2 \,dx \nonumber\\
			&\quad+2\int_{\Gamma_R}(x_2-c)\nu_2|\partial_\nu v|^2\,ds-2\int_{\Gamma}(x_2-c)\nu_2|\partial_{\nu}v|^2\,ds.  \label{eqn:rellich1}
		\end{align}
		Furthermore, integrating by parts, we have 
		\begin{align}
			&\quad  \int_{\Omega_R}(x_2-c)\partial_2|\nabla v|^2\,dx  \nonumber \\
			&= \int_{\partial\Omega_R}(x_2-c)\nu_2|\nabla v|^2\,ds-\int_{\Omega_R}|\nabla v|^2\,dx \nonumber \\
			&= \int_{\Gamma_R}(x_2-c)\nu_2|\nabla v|^2\,ds-\int_{\Gamma}(x_2-c)\nu_2|\nabla v|^2\,ds-\int_{\Omega_R}|\nabla v|^2\,dx. \label{eqn:rellich2}
		\end{align}
		We also obtain
		\begin{align}
			& \quad \int_{\Omega_R}k^2(x_2-c)\partial_2|v|^2 \,dx \nonumber\\
			&=\int_{\partial\Omega_R}k^2|v|^2(x_2-c)\nu_2\,ds-\int_{\Omega_R}k^2|v|^2\,dx \nonumber\\
			&= \int_{\Gamma_R}k^2|v|^2\nu_2(x_2-c)\,ds-\int_{\Gamma}k^2|v|^2\nu_2(x_2-c)\,ds-\int_{\Omega_R}k^2|v|^2\,dx.  \label{eqn:rellich3}
		\end{align}
		Substituting (\ref{eqn:rellich2}), (\ref{eqn:rellich3}) into (\ref{eqn:rellich1}) and using $$ -|\nabla v|^2+k^2|v|^2+2|\partial_2 v|^2=|\partial_2 v|^2-|\partial_1 v|^2+k^2|v|^2, $$ we have
		
		\begin{align*}
			&\quad  2\mbox{Re}\int_{\Omega_R}(x_2-c)\partial_2 \Bar{v} (\Delta v+k^2v)\,dx  \nonumber\\
			&= -2\int_{\Omega_R }|\partial_2 v|^2\,dx -\int_{\Gamma_R}(x_2-c)\nu_2|\nabla v|^2\,ds+\int_{\Gamma}(x_2-c)\nu_2|\nabla v|^2\,ds+\int_{\Omega_R}|\nabla v|^2\,dx \nonumber\\
			&\quad +\int_{\Gamma_R}k^2|v|^2\nu_2 (x_2-c)\,ds-\int_{\Gamma}k^2|v|^2\nu_2 (x_2-c)\,ds-\int_{\Omega_R}k^2|v|^2\,dx \nonumber\\
			&\quad +2\int_{\Gamma_R}(x_2-c)\nu_2|\partial_\nu v|^2\,ds-2\int_{\Gamma}(x_2-c)\nu_2|\partial_{\nu}v|^2\,ds \nonumber \\
			&=\left(\int_{\Gamma_R}-\int_{\Gamma}\right)(x_2-c)[-\nu_2|\nabla v|^2+\nu_2k^2|v|^2+2\mbox{Re}\,(\partial_2\bar{v}\partial_\nu v)]\,ds\nonumber\\
			&\quad+\int_{\Omega_R}|\nabla v|^2-k^2|v|^2-2|\partial_2 v|^2\,dx
		\end{align*}
		which completes the proof of Lemma \ref{eqn:rellichNohelmholtz}.
	\end{proof}
	\begin{corollary}\label{RemarkrellichNohelmholtz}
		If $v \in H_{\alpha}^2(\Omega_R)$ and $ v=0 \;\text{on} \;\Gamma$, then we have the Rellich identity 
		\begin{align*}
			\quad& 2{\rm{Re}}\,\int_{\Omega_R}(x_2-c)\partial_2 \Bar{v} (\Delta v+k^2v)\,dx+\int_{\Gamma}(x_2-c)\nu_2\left|\frac{\partial v}{\partial \nu}\right|^2\,ds+2\int_{\Omega_R}\left|\frac{\partial v}{\partial x_2}\right|^2\,dx\nonumber\\
			&=(R-c)\int_{\Gamma_R}\left|\frac{\partial v}{\partial x_2}\right|^2-\left|\frac{\partial v}{\partial x_1}\right|^2+k^2|v|^2\,ds+\int_{\Omega_R}|\nabla v|^2-k^2|v|^2\,dx,
		\end{align*}
		where $ c $ is a constant and $ \nu=(\nu_1,\nu_2) $ is the normal direction pointing to $ \Omega_{R} $.
	\end{corollary}
	\textbf{Proof of Lemma \ref{eqn:estimator1}.}  
	Define the half-plane 
	$$ U_{f_-}=\{x\in \mathbb{R}^2:x_2>f_-,\;x_1\in\mathbb{R}\}. $$
	For 
	$$v\in (C_0^\infty(\tilde{\Omega})\cap X_R) \subset (C_0^\infty( U_{f_-})\cap X_R),$$
	we have $$v(x)=\sum_{n \in \mathbb{Z}}v_n(x_2)e^{i\alpha_n x_1},\quad x_2>\Gamma_{\max}, $$ 
	where $v_n(x_2)  $ is the Fourier coefficient of $ v(x) $.
	Thus,
	\begin{align}
		|v_n(R)|^2 
		&= \int_{f_-}^{R}\frac{\partial}{\partial x_2}|v_n(x_2)|^2\,dx_2 \nonumber\\
		&=\int_{f_-}^{R}\frac{\partial}{\partial x_2}( v_n(x_2)\Bar{v}_n(x_2))\,dx_2 \nonumber\\
		&=2\mbox{Re}\,\int_{f_-}^{R}\Bar{v}_n(x_2)\frac{\partial}{\partial x_2}v_n(x_2)\,dx_2. \label{eqn:unR}
	\end{align}
	Using Cauchy-Schwarz inequality and $ (\ref{eqn:unR}) $, we have
	\begin{align}
		\|v\|_{H_{\alpha}^{1/2}(\Gamma_R)}^2& = \sum_{n \in \mathbb{Z}}(k^2+\alpha_n^2)^{1/2}|v_n(R)|^2 \nonumber\\
		&\leq 2\sum_{n \in \mathbb{Z}}(k^2+\alpha_n^2)^{1/2}\int_{f_-}^{R}|v_n(x_2)|\left|\frac{\partial}{\partial x_2}v_n(x_2)\right|\,dx_2 \nonumber\\
		&\leq 2\sum_{n \in \mathbb{Z}}(k^2+\alpha_n^2)^{1/2}\left(\int_{f_-}^{R}|v_n(x_2)|^2\,dx_2\right)^{1/2}\left(\int_{f_-}^{R}\left|\frac{\partial}{\partial x_2}v_n(x_2)\right|^2\,dx_2 \right)^{1/2}\nonumber\\
		&\leq 2\left(\sum_{n \in \mathbb{Z}}(k^2+\alpha_n^2)\int_{f_-}^{R}|v_n(x_2)|^2\,dx_2\right)^{1/2}\left(\sum_{n \in \mathbb{Z}}\int_{f_-}^{R}\left|\frac{\partial}{\partial x_2}v_n(x_2)\right|^2\,dx_2\right)^{1/2} \nonumber\\	
		& \leq \sum_{n \in \mathbb{Z}}(k^2+\alpha_n^2)\int_{f_-}^{R}|v_n(x_2)|^2\,dx_2+\sum_{n \in \mathbb{Z}}\int_{f_-}^{R}\left|\frac{\partial}{\partial x_2}v_n(x_2)\right|^2\,dx_2. \label{eqn:NormHVsX}
	\end{align}
	Note that in the last inequality we have used inequality $a+b\geq2\sqrt{ab}$ for $a,b>0$.
	Then we begin to prove the right-hand side of (\ref{eqn:NormHVsX}) is exactly $ \|v\|_X^2 /(2\pi)$.
	Define $ U_{f_-}^R=(0,2\pi)\times (f_-,R) $. 
	It follows that
	\begin{align}
		\|v\|_{L^2\left( \Omega_R\right)}^2&=\left\|\sum_{n \in \mathbb{Z}}v_n(x_2) e^{i\alpha_n x_1}\right\|_{L^2\left(\Omega_R\right)}^2 \nonumber\\
		&=\int_{U_{f_-}^R}\left|\sum_{n \in \mathbb{Z}}v_n(x_2) e^{i\alpha_n x_1}\right|^2\,dx \nonumber\\
		&=\int_{f_-}^{R}\int_{0}^{2\pi}\left(\sum_{n \in \mathbb{Z}}v_n(x_2) e^{i\alpha_n x_1}\right)\left(\sum_{m \in \mathbb{Z}}\bar{v}_m (x_2) e^{-i\alpha_m x_1}\right)\,dx_1dx_2 \nonumber\\
		&=2\pi\sum_{n \in \mathbb{Z}} \int_{f_-}^{R}|v_n(x_2)|^2\,dx_2. \label{l2}
	\end{align}
	Besides, similar to the above calculation, we have
	\begin{align}
		\|\nabla v\|_{L^2\left( \Omega_R\right)}^2&=\left\|\sum_{n \in \mathbb{Z}}v_n(x_2) i\alpha_ne^{i\alpha_n x_1}\right\|_{L^2\left(\Omega_R\right)}^2 +\left\|\sum_{n \in \mathbb{Z}}\frac{\partial}{\partial x_2}v_n(x_2) e^{i\alpha_n x_1}\right\|_{L^2\left(\Omega_R\right)}^2 \nonumber\\
		&=\int_{U_{f_-}^R}\left|\sum_{n \in \mathbb{Z}} \alpha_n v_n(x_2) e^{i\alpha_n x_1}\right|^2\,dx+\int_{U_{f_-}^R}\left|\sum_{n \in \mathbb{Z}}\frac{\partial}{\partial x_2}v_n(x_2) e^{i\alpha_n x_1}\right|^2\,dx\nonumber \\
		&=2\pi\sum_{n \in \mathbb{Z}}\alpha_n^2\int_{f_-}^{R}|v_n(x_2)|^2\,dx_2+2\pi\sum_{n \in \mathbb{Z}}\int_{f_-}^{R}\left|\frac{\partial}{\partial x_2}v_n(x_2)\right|^2\,dx_2.\label{nabl2}
	\end{align}
	Inserting (\ref{l2}) and (\ref{nabl2}) into (\ref{eqn:NormHVsX}) and using the definition of $ \|v\|_{X_R}$, we have
	\begin{eqnarray*}
		\|v\|_{H_{\alpha}^{1/2}(\Gamma_R)}\leq\|v\|_{X_R}/(2\pi)^{1/2}.   
	\end{eqnarray*}
	Since set $ \{v|_{\Omega_{R}}:v\in C_0^\infty(\tilde{\Omega})\cap X_R\} $ is dense in $ X_R $, we complete the proof.
	$\hfill\Box$
	
	\textbf{Proof of Lemma \ref{eqn:estimator2}.}
	For $ x\in \Omega_R $, we have
	\ben
	v(x)=\int_{f(x_1)}^{x_2}\frac{\partial v(x_1,\tau)}{\partial \tau} \,d\tau +v(x_1,f(x_1)).
	\enn
	Therefore, by Cauchy-Schwarz inequality we have
	\begin{align*}
		|v(x_1,x_2)|^2 &\leq 2\left|\int_{f(x_1)}^{x_2}\frac{\partial v(x_1,\tau)}{\partial \tau} \,d\tau \right|^2+2|v(x_1,f(x_1))|^2\\
		&\leq 2(x_2-f(x_1))\int_{f(x_1)}^{R}\left|\frac{\partial v(x_1,\tau)}{\partial \tau}\right|^2 \,d\tau+2|v(x_1,f(x_1))|^2.
	\end{align*}
	Thus, 
	\begin{align}
		\|v\|^2_{L^2(\Omega_R)}&=\int_{\Omega_R}|v|^2\,dx=\int_{0}^{2\pi}\int_{f(x_1)}^{R}|v(x_1,x_2)|^2\,dx_2 dx_1 \nonumber\\
		&\leq \int_{0}^{2\pi}\left[\left(\int_{f(x_1)}^{R}2(x_2-f(x_1))\,dx_2\right)\left(\int_{f(x_1)}^{R}\left|\frac{\partial v(x_1,\tau)}{\partial \tau}\right|^2 \,d\tau \right)\right]\,dx_1 \nonumber\\
		&\quad+\int_{0}^{2\pi} 2(R-f_-)|v(x_1,f(x_1))|^2\,dx_1 \label{lem4l2}
	\end{align}
	The first term on the right-hand side of the above formula can be estimated by
	\begin{align}
		&\quad\int_{0}^{2\pi}\left[\left(\int_{f(x_1)}^{R}2(x_2-f(x_1))\,dx_2\right)\left(\int_{f(x_1)}^{R}\left|\frac{\partial v(x_1,\tau)}{\partial \tau}\right|^2 \,d\tau \right)\right]\,dx_1\nonumber\\
		&\leq \int_{0}^{2\pi}(R-f(x_1))^2\left(\int_{f(x_1)}^{R}\left|\frac{\partial v(x_1,\tau)}{\partial \tau}\right|^2 \,d\tau\right) \,dx_1\nonumber\\
		&\leq (R-f_-)^2 \int_{0}^{2\pi} \int_{f(x_1)}^{R}\left|\frac{\partial v(x_1,\tau)}{\partial \tau}\right|^2 \,d\tau \,dx_1\nonumber\\
		&=  (R-f_{-})^2\left\|\frac{\partial v}{\partial x_2}\right\|^2_{L^2(\Omega_R)}. \label{lem41}
	\end{align}
	On the other hand,
	\begin{align}
		\|v\|^2_{L^2(\Gamma)}&=\int_{\Gamma}|v(x_1,x_2)|^2\;ds(x)\nonumber\\
		&=\int_{0}^{2\pi}|v(x_1,f(x_1))|^2\sqrt{1+(f'(x_1))^2}\;dx_1\nonumber\\
		&\geq \int_{0}^{2\pi}|v(x_1,f(x_1))|^2\;dx_1 \label{lem42}
	\end{align}
	Inserting (\ref{lem41}) and (\ref{lem42}) into (\ref{lem4l2}), we complete the proof.
	$\hfill\Box$
	
	\bibliographystyle{siam}
	\bibliography{refs}
\end{document}